\documentclass[10.5pt]{article}
\usepackage{enumerate,mathrsfs}
\usepackage{booktabs}
\usepackage{epic,epsf,epsfig,graphicx,color}
\usepackage{amsmath}
\usepackage{amssymb}
\usepackage{amsbsy}
\usepackage{amsfonts}
\usepackage{indentfirst}
\usepackage{multirow}
\newtheorem{theorem}{Theorem}[section]%
\newtheorem{lem}[theorem]{Lemma}%
\newtheorem{cor}[theorem]{Corollary}%
\newtheorem{prop}[theorem]{Proposition}%
\newtheorem{prob}[theorem]{Problem}%
\newtheorem{defi}[theorem]{Definition}%

\newtheorem{con}[theorem]{Conjecture}

\newtheorem{example}{Example}[section]
\newenvironment{exam}{\begin{example} \rm}{\end{example}\smallskip}

  \def\G{\Gamma}

\def\f{\noindent}

\def\Cay{\hbox{\rm Cay}}

\newcommand{\qed}{\mbox{\raisebox{0.7ex}{\fbox{}}} \vspace{4truemm}}

\def\mz{{\mathbb Z}}

\begin{document}
\title{A perfect matching reciprocity method for embedding multiple hypercubes in an augmented cube: Applications to Hamiltonian decomposition and fault-tolerant Hamiltonicity}
\footnotetext[1]{dwyang@bupt.edu.cn(D.-W. Yang), zhanghybupt@126.com(H. Zhang), rxhao@bjtu.edu.cn(R.-W. Hao), hsiehsy@mail.ncku.edu.tw(S.-Y. Hsieh)}

\author{ \vspace{1em} 
Da-Wei Yang$^{\rm a,b}$, Hongyang Zhang$^{\rm a}$, Rong-Xia Hao$^{\rm c}$, Sun-Yuan Hsieh$^{\rm d}$\\
\vspace{0.5em} 
{\small\em
\begin{tabular}{@{}l@{}}
    $^{\rm a}$School of Mathematical Sciences, Beijing University of Posts and Telecommunications, Beijing, 100876, P.R. China\\
    $^{\rm b}$Key Laboratory of Mathematics and Information Networks (BUPT), Ministry of Education, Beijing, 100876, P.R. China\\
    $^{\rm c}$School of Mathematics and Statistics, Beijing Jiaotong University, Beijing 100044, P.R. China.\\
    $^{\rm d}$Department of Computer Science and Information Engineering, National Cheng Kung University, Tainan 70101, Taiwan.
\end{tabular}
}
}
\date{}
\maketitle

\begin{abstract}
This paper focuses on the embeddability of hypercubes in an important class of Cayley graphs, known as augmented cubes. An $n$-dimensional augmented cube $AQ_n$ is constructed by augmenting the $n$-dimensional hypercube $Q_n$ with additional edges, thus making $Q_n$ a spanning subgraph of $AQ_n$.
Dong and Wang (2019) first posed the problem of determining the number of $Q_n$-isomorphic subgraphs in $AQ_n$, which still remains open. By exploiting the Cayley properties of $AQ_n$, we establish a lower bound for this number. What's more, we develop a method for constructing pairs of $Q_n$-isomorphic subgraphs in $AQ_n$ with the minimum number of common edges. This is accomplished through the use of reciprocal perfect matchings, a technique that also relies on the Cayley property of $AQ_n$.
As an application, we prove that $AQ_n$ admits  $n-1$ edge-disjoint Hamiltonian cycles when $n\geq3$ is odd and $n-2$ cycles when $n$ is even, thereby confirming a conjecture by Hung (2015) for the odd case. Additionally, we prove that $AQ_n$ has a fault-free cycle of every even length from $4$ to $2^n$ with up to $4n-8$ faulty edges, when each vertex is incident to at least two fault-free edges. This result not only provides an alternative proof for the fault-tolerant Hamiltonicity of established by Hsieh and Cian (2010), but also extends their work by demonstrating the fault-tolerant bipancyclicity of $AQ_n$.

\bigskip
\noindent{\bf Keywords:} Augmented cube; Hypercube; Cayley graphs; Perfect Matching, Hamiltonicity

\end{abstract}

\section{Introduction}

This paper focuses on a family of Cayley graphs known as augmented cubes, which play a significant role in large-scale interconnection networks and computer science. In multiprocessor systems, an interconnection network serves as the underlying topology, and it is usually represented as a graph in which vertices and edges correspond to nodes and communication links, respectively~\cite{DFH}. Throughout this paper, we adopt standard graph-theoretic terminology and notation~\cite{BM}.

Among the various interconnection networks proposed in the literature, the hypercube $Q_n$
stands out as one of the most efficient topologies due to its desirable properties, including recursive construction, low diameter, high connectivity, and symmetry~\cite{SS}.
To further enhance performance, particularly in reducing the diameter while preserving the hypercube's advantageous features, numerous variant networks were put forward successively~\cite{CS,ZFJZ}.
The augmented cube $AQ_n$ is one of these variants, which achieves a diameter roughly half that of the traditional hypercube $Q_n$~\cite{CS}. Beyond its improved diameter, the augmented cube retains many of the hypercube's favorable characteristics while offering superior embeddability--a property not shared by other hypercube variants~\cite{MKW,ZZC}. These advantages have made $AQ_n$ a subject of considerable interest in interconnection network research, prompting numerous studies and leaving several open problems unresolved. In the following, we summarize several related open problems, which motivate our selection of $AQ_n$ as the primary research subject of this work.

\subsection{Hypercubes in an augmented cube}

Proposed by Choudum and Sunitha~\cite{CS}, the augmented cube $AQ_n$ is constructed by augmenting the $n$-dimensional hypercube
with additional edges. These augmentation not only improve the network's robustness but also provide a richer structure for various algorithmic applications, including routing, broadcasting, and embedding~\cite{CS,WSQ}. Since a hypercube $Q_n$ is a spanning subgraph of $AQ_n$, the augmented cube $AQ_n$ inherits many favorable properties of the hypercube, including maximum connectivity, optimal wide diameter, and support for time-efficient routing and broadcasting algorithms with linear complexity~\cite{CS}. In particular, an augmented cube $AQ_n$ has high symmetric properties as a hypercube $Q_n$; both of them are Cayley graphs on elementary abelian 2-groups~\cite{CS2}. These compelling characteristics have spurred extensive research into the topological properties of augmented cubes, such as reliability~\cite{CLQS}, fault tolerance~\cite{ML,ZXWZY}, diagnosability~\cite{GWXLL},
and embeddability~\cite{MKW, WSQ}, making it a promising candidate for large-scale parallel and distributed systems~\cite{ZZC}.

Embeddability is a critical issue in the design of an interconnection network, as it determines the capability of a network to incorporate other existing network topologies~\cite{X}.
Choudum and Sunita~\cite{CS} demonstrated that the augmented cube $AQ_n$ possesses several remarkable embedding properties not found in hypercubes or their variants. Specifically, they proved that an augmented cube of
dimension $n$ contains two edge-disjoint complete binary trees on $2^n-1$ vertices both rooted at the same vertex, two edge-disjoint spanning binomial trees, and cycles of all possible length from 3 to $2^n$. Other better embeddability of augmented cubes has also been presented in some research, such as~\cite{MKW,MG}.

Cycles serve as fundamental networks for parallel and distributed computation, being suitable for the development of simple algorithms with low communication cost~\cite{LHJ}. Embedding cycles into augmented cubes is vital and has been investigated in the literature, see~\cite{QM,WSQ} for example. Specially, Dong and Wang~\cite{DW} have studied the cycles in an augmented cube of the smallest possible length, and enumerated the numbers of cycles of length 3 or 4. Their analytical approach leverages the critical structural relationship between hypercubes and augmented cubes, that is the existence of $Q_n$ as a spanning subgraph of $AQ_n$. Particularly, they posed the following problem.
\begin{prob} {\rm (\cite{DW})}
How many $n$-dimensional hypercubes in an $n$-dimensional augmented cube?
\end{prob}

Let $n$ be an integer at least 2, and let
$$f(n)=1+\sum\limits_{1\leq k \leq n-1}\
 \sum\limits_{2\leq j_1 < j_2<\cdots <j_k\leq n}{(j_k-j_{k-1})\cdots (j_2-j_1)j_1}.$$
Through using the Cayley properties of the hypercube and the augmented cube,
we prove that there are at least $f(n)$ distinct $n$-dimensional hypercubes in $AQ_n$.
This gives a lower bound of the number of subgraphs in $AQ_n$ that are isomorphic to $Q_n$.

The function $f(n)$ quantifies the abundance of $Q_n$-isomorphic subgraphs in augmented cubes $AQ_n$. While $AQ_n$ contains numerous such subgraphs, many exhibit only trivial structural variations. For most applications, selecting between two nearly identical subgraphs offers no substantive advantage over using a single subgraph. A more valuable configuration involves identifying $Q_n$-isomorphic subgraphs in $AQ_n$ with minimal edge overlap.
Given that $AQ_n$ is a $(2n-1)$-regular graph while $Q_n$ is $n$-regular, we observe that when two $Q_n$-isomorphic subgraphs in $AQ_n$ minimize their edge intersection, this intersection forms a perfect matching (see Section 3.2). Building on this foundation, we develop a systematic construction of
$Q_n$-isomorphic subgraphs with minimal edge overlap by exploiting the Cayley properties of
$AQ_n$. The core of our approach involves a novel reciprocal perfect matching method that generates such distinct subgraph pairs. The remainder of this paper demonstrates how this construction provides crucial insights for solving several fundamental problems concerning $AQ_n$.

\subsection{Hamiltonian cycles in an augmented cube}

In contrast to determining the shortest cycle in a network, researchers tend to focus more on investigating the longest cycles present within the network. A classic problem in this context is the Hamiltonicity of a graph. In a graph $G$, a Hamiltonian cycle is a cycle that contains every vertex of $G$, and $G$ is called Hamiltonian if it contains a Hamiltonian cycle.
It is known that both hypercubes and augmented cubes are Hamiltonian graphs~\cite{CS}.
In a Hamiltonian graph, there are usually many Hamiltonian cycles.
If two Hamiltonian cycles in a graph have no common edge, then they are called {\em edge-disjoint}.
Edge-disjoint Hamiltonian cycles are often abbreviated to EDHCs.
EDHCs have many applications and they can provide an advantage for algorithms using ring structures~\cite{Hung2}. Furthermore, it can also be used to solve the problem of all-to-all communication algorithms and fault-tolerant routing~\cite{HBA,PWPC}. These applications motivated researchers to explore EDHCs in a network. For example, Hung, L\" v and et al.~\cite{Hung2,Hung3,LW,YXFL} proved that there exist two EDHCs in many variants of hypercubes,
including crossed cubes, locally twisted cubes, spined cubes and balanced hypercubes. Researchers have persistently pursued the identification of more EDHCs within networks. For instance, the existence of three EDHCs in crossed cubes and locally twisted cubes has been demonstrated in~\cite{PWPC,Pai}.

Clearly, an $n$-regular graph $G$ has at most $\lfloor \frac{n}{2}\rfloor$ EDHCs.
If the number of EDHCs in $G$ achieves this theoretical upper bound, then these Hamiltonian cycles are said to form a {\em Hamiltonian decomposition} of $G$. Bae and Bose~\cite{BB} have proved that an $n$-dimensional hypercube admits a Hamiltonian decomposition.
In 2015, Hung~\cite{Hung3} conjectured that this also holds for augmented cubes.
\begin{con}\label{con=1}{\rm (\cite{Hung3})}
An $n$-dimensional augmented cube admits a Hamiltonian decomposition.
\end{con}
As our best knowledge, this conjecture is still open.
After proving that there are many hypercubes in an augmented cube,
we select two $n$-dimensional hypercubes in an augmented cube $AQ_n$ whose
edges sets intersect in a perfect matching. Through using the Hamiltonian decompositions
of these two hypercubes, we prove that there are $m$ EDHCs in an $AQ_n$, where
$m=n-1$ if $n$ is odd and $m=n-2$ otherwise. This not only extends the previous best result on EDHCs of $AQ_n$ presented by Hung~\cite{Hung1} who constructed two EDHCs, and also confirms Conjecture~\ref{con=1}
when $n$ is odd.

\subsection{Edge-fault-tolerant Hamiltonicity of augmented cubes}

Element (edge or vertex) failure is inevitable when an interconnection network is put in use. Therefore, the fault-tolerant capacity of a network is a critical issue. A graph or network is said to be {\em faulty} if it has at
least one faulty vertex or edge.
An important research focus involves investigating the structural properties of the remaining functional subgraphs after moving all faulty elements, such as the fault-tolerant cycle embedding problem in networks~\cite{XM}.
We focus on the scenario where only edge faults exist in a graph, with the research typically conducted under two distinct modeling frameworks~\cite{HC}. One is the random fault model: faults might happen anywhere in a graph without any restriction~\cite{HHC}. The other one is the conditional fault model: each vertex is incident to at least two or more fault-free edges~\cite{HC}. The conditional fault model is practically motivated that in real-world network applications, all edges incident to the same vertex fail at the same time is almost impossible~\cite{Harary}.

The edge-fault-tolerant Hamiltonian problem investigates whether the remaining subgraph maintains Hamiltonian properties after the removal of certain faulty edges from the original graph. A cycle in a faulty graph is {\em fault-free} if it contains no faulty edges. Chan and Lee~\cite{CL} first studied the fault-tolerant Hamiltonicity of a hypercube with a set $F$ of faulty edges. They showed that there
exists a fault-free Hamiltonian cycle in $Q_n$
with $n\geq 3$ and $|F|\leq 2n-5$, under the
constraint that each vertex is incident to at least two fault-free edges. Later, Liu and Wang~\cite{LW} extended Chan and Lee's work~\cite{CL}, by proving that under one more constraint, a $Q_n$ with $n\geq 5$ also remains Hamiltonicity even if $|F|\leq 3n-8$ (see Section~5 for detail).

The edge-fault-tolerant Hamiltonicity of an augmented cube has also been investigated.
Let $F$ be the set of faulty edges in $AQ_n$.
Under the conditional fault model, which assumes that each vertex is incident to at least two fault-free edges,
M\v ekuta and Gregor~\cite{MG} established the existence of a fault-free Hamiltonian cycle
in $AQ_n$ when $n\geq4$ and $|F|\leq 3n-7$.
Additionally, they provided a secondary result-though less refined-stating that $AQ_n$ contains a fault-free Hamiltonian cycle when $|F|\leq 4n-10$ and each vertex is incident to at least $n+1$ fault-free edges. Notably, they posed a problem as follows.

\begin{prob}{\rm (\cite{MG})}
If $|F|\leq 4n-10$ and each vertex is incident to at least $2$ fault-free edges, does $AQ_n$ contain a fault-free Hamiltonian cycle?
\end{prob}

The reason we consider this result less refined is that $4n-10$ is not an optimal upper bound. Hsieh and Cian~\cite{HC2} gave an optimal upper bound under the conditional fault model. They proved that $AQ_n$ contain a fault-free Hamiltonian cycle when
$|F|\leq 4n-8$ and each vertex is incident to at least two fault-free edges. An counterexample of $AQ_n$ with $4n-7$ faulty edges was also constructed. Their proof employs a mathematical induction, leveraging the recursive structure of $AQ_n$.
In this paper, we present an alternative proof of Hsieh and Cian's result.
Rather than induction, our approach utilizes perfect matching reciprocity to construct $Q_n$-isomorphic subgraphs within $AQ_n$ that contains as few faulty edges as possible, offering a new perspective on the problem.
Moreover, we extend their result by showing that
$AQ_n-F$ not only remains Hamiltonian but also contains cycles of every even length from 4 to
$2^n$ under the same fault constraints. This result connects to the broader study of fault-tolerant (bi)pancyclicity in interconnection networks; for related work, we refer readers to~\cite{CHF}.

\medskip

The remainder of this paper is organized as follows. Section 2 introduces some definitions and notations. In Section 3, leveraging the Cayley property of $AQ_n$, we systematically investigate how to identify subgraphs isomorphic to $Q_n$ within $AQ_n$. Special emphasis is placed on finding pairs of such subgraphs with the minimum number of common edges. Sections 4 and 5 then present two applications of this methodology: the construction of edge-disjoint Hamiltonian cycles in $AQ_n$ and the edge-fault-tolerant Hamiltonicity of $AQ_n$. Section 6 is the conclusion.

\section{Preliminaries}

To study Hamiltonian cycles and the subgraphs isomorphic to $Q_n$ in an augmented cube, our main idea is using the Cayley properties of hypercubes and augmented cubes. Therefore, we need a few notations of elementary group theory. To avoid ambiguity
in mathematical symbol, we use $\Gamma$ to denote a group and $G$ to represent a graph.
Unless stated otherwise, we follow~\cite{BM} and~\cite{KS} for terminology and definitions related to graphs and groups, respectively.

\subsection{Fundamental group and graph terminologies}

All groups in this paper are finite. Let $\G$ be a group with identity $e$. The group $\G$ is {\em abelian} if $ab=ba$ for any $a,b\in \G$. For each $s\in \G$,
we use $s^{-1}$ to denote the inverse of $s$, i.e., $ss^{-1}=s^{-1}s=e$.
For a positive integer $n$ and a prime $p$, we use $\mz_p^n$ to denote an elementary abelian $p$-group of order $p^n$.

Let $\G$ be an abelian group. For a nonempty subset $S=\{x_1,\ldots,x_n\}$ of $\G$,
we use $\langle S\rangle$ to denote the subgroup generated by $S$,
that is,
$$\langle S\rangle=\{x_1^{z_1}\cdots x_n^{z_i}~|~x_i\in S, z_i\in \mz, 1\leq i\leq n\},$$
where $\mz$ is the set of integers. The subset $S$ is said to be a {\em generating subset} of $\langle S\rangle$. We also write $\langle S\rangle=\langle x_1,\ldots,x_n\rangle$.
Clearly, $\langle S\rangle\subseteq \Gamma$ if and only if $x_i\in \Gamma$ for each $1\leq i\leq n$. A group may has many generating subsets. A generating subset $S$ of group $\G$
is {\em minimal} if $S\backslash\{x\}$ cannot generates $\G$ for any $x\in S$.
For example, in an elementary abelian $2$-group
$$\mz_2^n=\langle a_1,a_2,\ldots, a_n~|~a_i^2=e,a_ia_j=a_ja_i,1\leq i,j\leq n\rangle,$$
the sets $\{a_1,a_2,\ldots, a_n\}$ and $\{a_1a_n,a_2,\ldots,a_n\}$ are two minimal generating subsets of $\mz_2^n$.


\medskip

Throughout this paper, all graphs are finite, simple and undirected.
Let $G$ be a graph with vertex set $V(G)$ and edge set $E(G)$. We use $d_G(u)$ to denote the degree of the vertex $u$ in $G$, that is the number of edges incident with $u$. If $d_G(u)=k$ for each $u\in V(G)$, then $G$ is said to be a {\em $k$-regular graph}. Given a subset $F$ of $V(G)$, the subgraph of $G$ induced by $F$, denoted by $G[F]$, is the graph whose vertex set is $F$ and edge set is $\{(u, v)\in E(G)~|~u,v\in F\}$. The notation $G-F$ represents the subgraph of $G$ induced by $V(\G)\backslash F$.
For a subset $E$ of $E(G)$, the {\em edge-induced subgraph} $G[E]$ is the subgraph of $G$ whose edge set is $E$ and whose vertex set consists of all ends of the edges of $E$. A graph $H$ is a {\em spanning subgraph} of $G$ if $V(H)=V(G)$ and $E(H)\subseteq E(G)$.

Let $C=(u_{1},u_{2},\ldots,u_{s},u_1)$ be a cycle in a graph $G$. For an edge $(u_i,u_{i+1})$ in $C$,
we use $C-(u_i,u_{i+1})$ to denote the path $(u_{i+1},\ldots,u_s,u_1,\ldots,u_{i})$,
and use $(u_1,\ldots,u_i)+(u_i,u_{i+1})$ to denote the path $(u_1,\ldots,u_i,u_{i+1})$, with the understanding that this represents a cycle if
$u_{i+1}=u_1$.

Two graphs $G$ and $H$ are {\em disjoint} if they have no vertex in common, and {\em edge-disjoint} if they have no edge in common.
An {\em isomorphism} from $G$ to $H$ is a bijection $\phi: V(G)\rightarrow V(H)$ such that $(u, v)\in E(G)$ if and only if $(\phi(u),\phi(v))\in E(H)$. The graphs $G$ and $H$ are {\em isomorphic}, write $G\cong H$, if there is an isomorphism from $G$ to $H$.
An isomorphism from $G$ to itself is called an {\em automorphism} of $G$.
We say that $G$ is {\em vertex-transitive} if for any two vertices $u, v\in V(G)$, there exists an automorphism $\phi$ of $G$ such that $\phi(u)=v$.

The class of Cayley graphs plays an important role in studying vertex-transitive graphs.
Let $\G$ be a group with identity $e$. Given a subset $S$ of $\G$ such that
$S^{-1}=S$ and $e\notin S$, where $S^{-1}=\{s^{-1}~|~s\in S\}$,
the {\em Cayley graph} of $\G$ with respect to $S$, write $\Cay(\G, S)$, is the graph with vertex set $\G$ and edge set $\{(g, sg)~|~g\in \G, s\in S\}$. A Cayley graph $\Cay(\G, S)$ is connected if and only if $\langle S\rangle=\G$. Moreover, if $\langle S\rangle\neq \G$, then the disconnected graph $\Cay(\G, S)$ has $|\G|/|\langle S\rangle|$ components, in which each component is isomorphic to $\Cay(\langle S\rangle, S)$. It is well-known that a Cayley graph is vertex-transitive, while the converse is not true.

\begin{figure}
	\centering	\includegraphics[width=0.6\linewidth]{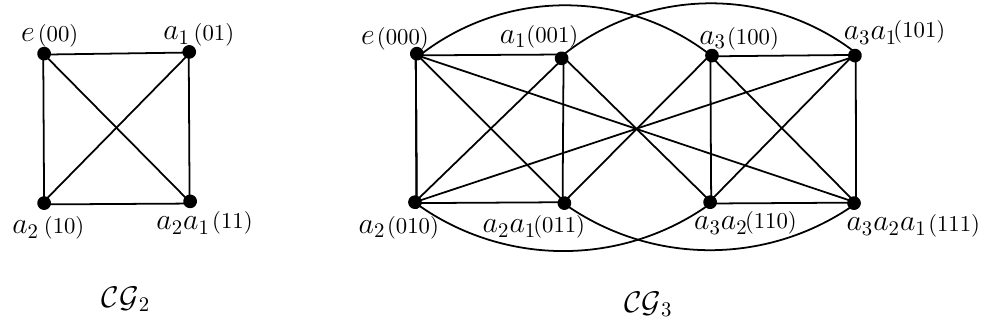}
	\caption{The Cayley graphs $\mathcal{CG}_{2}$ and $\mathcal{CG}_{3}$}
	\label{fig=0}
\end{figure}

\begin{exam}\label{exam=1}
Let $\mz_2^n=\langle a_1,a_2,\ldots, a_n\rangle$ be an elementary abelian $2$-group with $n\geq2$. Define a Cayley graph $\mathcal{CG}_{n}=\Cay(\mz_2^n,S)$ with
\begin{equation}\label{label=s}
S=\{a_1,a_2,\ldots,a_n, a_2a_1, a_3a_2a_1, \ldots , a_na_{n-1}\cdots a_1 \}.
\end{equation}
Then $\mathcal{CG}_{n}$ is a $(2n-1)$-regular graph. Since $\langle S\rangle=\mz_2^n$,
$\mathcal{CG}_{n}$ is connected. The graphs $\mathcal{CG}_{2}$ and $\mathcal{CG}_{3}$ are depicted in Figure~\ref{fig=0}.
\end{exam}

\subsection{Hypercubes and augmented cubes}

Let $n$ be a positive integer. An {\em $n$-dimensional hypercube $Q_n$} is a graph of order $2^{n}$. Each vertex is labeled with a binary string $x_{n}x_{n-1}\cdots x_{2}x_{1}$, and two vertices are adjacent if they have only one bit distinct.
The following proposition is well-known.

\begin{prop}\label{prop=2.3-1}
Let $n\geq2$. Then $Q_n\cong\Cay(\mz_2^n, T)$ for any  minimal generating subset $T$ of the elementary abelian $2$-group $\mz_2^n$.
\end{prop}

The \textit{$n$-dimensional augmented cube} $AQ_n$ was first defined in~\cite{CS}, which has the same vertex set as $Q_n$.

\begin{defi}\label{defi=1}
The $AQ_1$ is a complete graph $K_2$ of two vertices labeled with $0$ and $1$, respectively. For $n\geq2$, $AQ_n$ is obtained by taking two copies of $AQ_{n-1}$ with the label of each vertex preceded by $0$ $($or $1$, resp.$)$, denoted by $0AQ_{n-1}$ and $1AQ_{n-1}$, and adding $2\times 2^{n-1}$ edges between the two copies as follows. A vertex $u=0u_{n-1}u_{n-2}\cdots u_1$ of $0AQ_{n-1}$ is adjacent to a vertex $v=1v_{n-1}v_{n-2}\cdots v_1$ of $1AQ_{n-1}$ if and only if either of these conditions is satisfied.
\begin{enumerate}
\item [\rm (1)] $u_i=v_i$ for $1\leq i\leq n-1$; in this case, $(u, v)$ is called a hypercube edge of dimension $n$,
\item [\rm (2)] $u_i=1+v_i$ for $1\leq i\leq n-1$; in this case, $(u, v)$ is called an augmented edge of dimension $n$.
\end{enumerate}
\end{defi}

The augmented cubes $AQ_2$ and $AQ_3$ are indeed the graphs $\mathcal{CG}_2$ and $\mathcal{CG}_3$, respectively, as depicted in Figure~\ref{fig=0}.
For any vertex $u=x_n\cdots x_1$ in $AQ_n$, by Definition~\ref{defi=1}, its neighborhood is
\begin{equation}
\begin{split}\label{eq=1}
N_{AQ_n}(u)=\{&x_n\cdots x_2(x_1+1), x_n\cdots x_3(x_2+1)x_1,\ldots,(x_n+1)x_{n-1}\cdots x_2x_1,\\
&x_n\cdots x_3 (x_2+1) (x_1+1), x_n\cdots x_4(x_3+1)(x_2+1)(x_1+1),\ldots,\\
& (x_n+1)(x_{n-1}+1)\cdots (x_1+1)\}.
\end{split}
\end{equation}


Although an augmented cube has been proved to be a Cayley graph over an elementary abelian $2$-group, we still provide a proof to enhance the readability of the paper.

\begin{lem}\label{lem=1}
Let $n\geq2$. Then $AQ_n\cong\mathcal{CG}_n$.
\end{lem}

\f {\bf Proof.} Note that the vertex set of $\mathcal{CG}_n$ is the elementary abelian group
$\mz_2^n=\langle a_1,\ldots, a_n\rangle$, and each element of $\mz_2^n$ can be uniquely represented as $a^{x_n}_n\cdots a^{x_1}_1$, where $x_i\in \{0,1\}$ for $1\leq i\leq n$.
Define a map from $V(AQ_n)$ to $V(\mathcal{CG}_n)$ as follows:
$$\phi:u=x_n\cdots x_1 \mapsto a^{x_n}_n\cdots a^{x_1}_1,~ \forall u\in V(AQ_n).
$$
Then $\phi$ is a bijection from $V(AQ_n)$ to $V(\mathcal{CG}_n)$.
Since both $AQ_n$ and $\mathcal{CG}_n$ are $(2n-1)$-regular graphs, to show that $\phi$ is an isomorphism from $AQ_n$ to $\mathcal{CG}_n$, it suffices to show
$\phi[N_{AQ_n}(u)]= N_{\mathcal{CG}_n}(\phi(u))$ for any $u\in V(AQ_n)$.

Let $u=x_n\cdots x_1\in V(AQ_n)$. By Eq.~\eqref{eq=1}, we have
$$\begin{array}{ll}
\phi[N_{AQ_n}(u)]=\{&a_n^{x_n}\cdots a_2^{x_2}a_1^{x_1+1}, a_n^{x_n}\cdots a_3^{x_3}a_2^{x_2+1}a_1^{x_1},\ldots,
a_n^{x_n+1}a_{n-1}^{x_{n-1}}\cdots a_1^{x_1},a_n^{x_n}\cdots a_3^{x_3}a_2^{x_2+1}a_1^{x_1+1},\\
&a_n^{x_n}\cdots a_{4}^{x_{4}}a_3^{x_3+1}a_{2}^{x_{2}+1}a_1^{x_1+1},  \ldots, a_n^{x_n+1}a_{n-1}^{x_{n-1}+1}\cdots a_1^{x_1+1}\}.
\end{array}$$
On the other hand, since $\phi(u)=a_n^{x_n}\cdots a_1^{x_1}$, by Example~\ref{exam=1} we have
$\phi[N_{AQ_n}(u)]= N_{\mathcal{CG}_n}(\phi(u))$, as required.
\hfill\qed

The following two proposition come from the proof of~\cite[Theorem 3.2]{DW}.

\begin{prop}\label{prop=2.3-3}
Let $n\geq3$, and let $C$ be a $4$-cycle of $AQ_n$.
\begin{itemize}
\item [\rm (1)] If $C$ consists of exactly one augmented edge, then $C$ has an augmented edge of dimension $3$ and three hypercubes edges of dimensions $1$, $2$, $3$, respectively.
\item [\rm (2)] If $C$ consists of exactly two hypercube edges and one of which is of dimension $1$, then the other hypercube edge is also of dimension $1$.
\item [\rm (3)] If $C$ consists of exactly two augmented edges of dimensions $j_1$ and $j_2$ with $j_1<j_2$ that are adjacent, then $j_2=j_1+2$ and the two hypercube edges have dimensions $j_2$ and $j_1+1$, respectively.
\item [\rm (4)] The $C$ cannot have exactly three augmented edges.
\end{itemize}
\end{prop}

\begin{prop}\label{prop=2.3-2}
Let $n\geq3$. There are $2^{n-2}(2n^2+5n-11)$ cycles of length $4$ in $AQ_n$ in total.
For each edge $e\in E(AQ_n)$, there are at most $2n+8$ cycles of length $4$ going through $e$
in $AQ_n$.
\end{prop}

\section{Hypercubes in an augmented cube}

In this section, we aim to find different subgraphs of $AQ_n$ that are isomorphic to $Q_n$. For convenience, we call these subgraphs as $Q_n$-isomorphic subgraphs.

\subsection{A lower bound of the number of $Q_n$-isomorphic subgraphs in $AQ_n$}

In a Cayley graph $G=\Cay(\G,S)$, for any subset $T\subseteq S$ satisfying $T^{-1}=T$, the edge-induced subgraph $G[E]$
with $E=\{(g,tg)~|~t\in T\}$ is also a Cayley graph of $G$ which is indeed $\Cay(\G,T)$. In light of Proposition~\ref{prop=2.3-1} and Lemma~\ref{lem=1}, both a hypercube $Q_n$ and an augmented cube $AQ_n$ can be represented as Cayley graphs over the elementary abelian $2$-group $\mz_2^n$. To identify subgraphs of $AQ_n$ that are isomorphic to $Q_n$, it suffices to locate Cayley subgraphs of $\mathcal{CG}_n=\Cay(\mz_2^n,S)$ that are isomorphic to $Q_n$. Given that $Q_n\cong\Cay(\mz_2^n,T)$ with $T$ a minimal generating subset of $\mz_2^n$,
one approach to finding such $Q_n$-isomorphic subgraphs within $\mathcal{CG}_n\cong AQ_n$ is to analyze the minimal generating subsets of $S$ (see Eq.~\eqref{label=s} for $S$).

\begin{theorem}\label{the=1}
In $AQ_n$ with $n\geq2$, there are at least $f(n)$ distinct $Q_n$-isomorphic subgraphs, where $$f(n)=1+\sum\limits_{1\leq k \leq n-1}\
 \sum\limits_{2\leq j_1 < j_2<\cdots <j_k\leq n}{(j_k-j_{k-1})\cdots (j_2-j_1)j_1}.$$
\end{theorem}

\f {\bf Proof.} We use Cayley graphs $\Cay(\mz_2^n,T)$
and $\Cay(\mz_2^n,S)$ to represent $Q_n$ and $AQ_n$, respectively, where $\mz_2^n=\langle a_1,\ldots,a_n\rangle$. The $T$ is a minimal generating subset of $\mz_2^n$ and $S=S_1\cup S_2$ with
$$S_1=\{a_1,a_2,\ldots,a_n\},~S_2=\{a_2a_1, a_3a_2a_1, \ldots , a_na_{n-1}\cdots a_1\}.$$
Our goal is to find all possible $T$s in $S$.
It is worth to note that a minimal generating subset of $\mz_2^n$ consists of exactly $n$ elements,
while not any $n$ elements of $\mz_2^n$ form a generating subset. The problem can be translated into finding $n$ elements $b_1,\ldots ,b_n$ in $S$ such that they can generated the group $\mz_2^n$, namely $\mz_2^n=\langle b_1,\ldots , b_n\rangle$.
Let $T=\{b_1,\ldots ,b_n\}$. We consider the following three cases depending on the intersection $T\cap S_2$.

\f {\bf Case~1:} $|T\cap S_2|=0$, i.e. $T\subseteq S_1$.

Since $|T|=|S_1|=n$, we have that $T=S_1=\{a_n,\ldots ,a_1\}$. Obviously, $S_1$ is a minimal generating subset of $\mz_2^n$. Now, we obtain a subgraph of $AQ_n$ that is isomorphic to $AQ_n$, which is the subgraph induced by the edge set
$$E=\{(g,sg)\,|\, g\in \mz_2^n,\ s\in S_1\}=\{(g, a_1 g),\ldots , (g, a_ng)\,|\,g\in \mz_2^n\}.$$

\f {\bf Case~2:} $|T\cap S_2|=1$, i.e. $T\cap S_2=\{a_ia_{i-1}\cdots a_1\}$ for some $i$ with $2\leq i\leq n$.

Since $T$ is a subset of $S$ of cardinality $n$, we have $T=\{a_1,a_2,\ldots,a_n,a_ia_{i-1}\cdots a_1\}\backslash \{a_j\}$ for some $1\leq j\leq n$,
and since $T$ is a generating subset of $\mz_2^n$,
we have $1\leq j\leq i$ (otherwise, $\mz_2^{n-1}\cong\langle T\rangle\neq \mz_2^n$).
Conversely, for each
$1\leq j\leq i$, the subset $\{a_1,a_2,\ldots,a_n,a_ia_{i-1}\cdots a_1\}\backslash \{a_j\}$ is clearly a minimal generating subset of $\mz_2^n$.
Hence we have
$$\begin{array}{ll}
T=&\{a_2,\ldots,a_{n}, a_ia_{i-1}\cdots a_1\}, \{a_1,a_3,\ldots, a_n, a_ia_{i-1}\cdots a_1\},\ldots,\\
&\{a_1,\ldots,a_{i-2}, a_i, a_{i+1},\ldots , a_n, a_{i} a_{i-1}\cdots a_1\},~{\rm or}\\
&\{a_1,\ldots,a_{i-1}, a_{i+1},\ldots , a_n, a_i a_{i-1}\cdots a_1\}.
\end{array}$$
In conclude, for each $2\leq i\leq n$, there exists $i$ distinct minimal generating subsets of $\mz_2^n$.
Hence this case yeilds $\sum\limits_{2\leq i\leq n}{i}$ distinct $Q_n$-isomorphic subgraphs.
\medskip

\f {\bf Case~3:} $|T\cap S_2|=k$ with $2\leq k\leq n-1$.

Assume that $T=\{a_1,\ldots ,a_n, a_{j_1}\cdots a_1, a_{{j_{2}}}\cdots a_1,\ldots ,a_{j_k} \cdots a_1\}\backslash \{a_{i_1},a_{i_2},\ldots,a_{i_k}\}$,
where $1\leq i_1 < i_2 < \cdots < i_k\leq n$
and $2\leq j_1 < j_2 < \cdots < j_k\leq n$.
Since $\langle T\rangle=\mz_2^n$, we have $i_k\leq j_k$.

Suppose $i_k\leq j_{k-1}$. Then $i_k<j_k$ because $j_{k-1}<j_k$. It forces that $\{a_{j_{k-1}+1},\ldots,a_{j_k},\ldots , a_{n}\} \subseteq T$. Moreover, $\{a_{j_{k-1}}\cdots a_1,a_{j_k}\cdots a_1\}\subseteq T$.
Since
$$a_{j_k}\cdots a_1 =a_{j_k}\cdot a_{j_k-1}\cdot \ \cdots \ \cdot (a_{j_{k-1}+1})\cdot (a_{j_{k-1}}\cdots a_1),$$
we have $$\langle a_{j_{k-1}+1},\ldots,a_{j_k},\ldots a_{n}, a_{j_{k-1}}\cdots a_1, a_{j_k}\cdots a_1 \rangle
=\langle a_{j_{k-1}+1},\ldots,a_{j_k},\ldots a_{n}, a_{j_{k-1}}\cdots a_1\rangle. $$
This means that $\langle T\rangle=\langle T\backslash \{a_{j_k}\cdots a_1\}\rangle$, which is contradict to
the minimality of $T$. Therefore, $j_{k-1}<i_k\leq j_k$.

Now, we turn to compare the integers $i_{k-1}$ and $j_{k-1}$. Suppose $i_{k-1}>j_{k-1}$. Then
$j_{k-1}<i_{k-1}<i_k\leq j_k$. It forces that for the $k$ elements in $T\cap S_1=\{a_{j_1}\cdots a_1, a_{{j_{2}}}\cdots a_1,\ldots ,a_{j_k} \cdots a_1\}$, both $a_{i_{k-1}}$ and
$a_{i_k}$ only appear in the expression of the element $a_{j_k} \cdots a_1$.
Hence,
$$a_{i_k}\notin \langle a_1,\ldots, a_{i_{k-1}-1},a_{i_{k-1}+1},\ldots, a_{i_k-1},a_{i_{k}+1},\ldots,a_n,a_{j_1}\cdots a_1, a_{{j_{2}}}\cdots a_1,\ldots ,a_{j_k} \cdots a_1\rangle.$$
Since $T\subseteq \{a_1,\ldots, a_{i_{k-1}-1},a_{i_{k-1}+1},\ldots, a_{i_k-1},a_{i_{k}+1},\ldots,a_n,a_{j_1}\cdots a_1, a_{{j_{2}}}\cdots a_1,\ldots ,a_{j_k} \cdots a_1\}$,
we have $a_{i_k}\notin \langle T\rangle$, which is a contradiction because $T$ is a generating subset of $\mz_2^n$. Hence $i_{k-1}\leq j_{k-1}$.

Following a method analogous to the one above, by sequentially examining two consecutive integers in the sequence $(i_{k-1}, j_{k-2}, i_{k_2},\ldots, i_2,j_1,i_1)$, we ultimately arrive at the following inequation
\begin{equation}\label{eq=2}
1\leq i_1\leq j_{1}<i_2\leq j_2<\cdots <i_{k-1}\leq j_{k-1}<i_{k}\leq j_k\leq n,~j_1\geq 2.
\end{equation}

Conversely, let $T=\{a_1,\ldots ,a_n, a_{j_1}\cdots a_1, a_{{j_{2}}}\cdots a_1,\ldots ,a_{j_k} \cdots a_1\}\backslash \{a_{i_1},a_{i_2},\ldots,a_{i_k}\}$
satisfying the inequation~\eqref{eq=2}.
Note that $1\leq i_1\leq j_{1}$ and $2\leq j_1$.
If $i_1<j_1$, then $\{a_1,\ldots , a_{i_1-1}, a_{i_1+1},\ldots, a_{j_1}, a_{j_1}\cdots a_1\}\subseteq T$ (if specially $i_1=1$ then $\{a_1,\ldots , a_{i_1-1}\}$ is regarded as an empty set),
and since
$a_{i_1}=(a_{j_1}\cdots a_1)\cdot a_1\cdot \ \cdots \ \cdot a_{i_1-1}\cdot a_{i_1+1}\cdot \ \cdots\ \cdot a_{j_1},$
we have $a_{i_1}\in \langle T\rangle$. If $i_1=j_1$,
then $\{a_1,\ldots , a_{i_1-1}, a_{i_1}\cdots a_1\}\subseteq T$
and so $a_{i_1}=(a_{i_1}\cdots a_1)\cdot a_1\cdot \ \cdots \ \cdot a_{i_1-1}\in \langle T\rangle$.

Similarly, for any integer $2\leq s\leq k$, we have $i_{s-1}\leq j_{s-1}<i_s\leq j_s$.
If $i_{s}=j_s$, then we have
$\{a_{i_{s-1}+1},\ldots,a_{j_{s-1}+1},\ldots,a_{i_s-1}, a_{i_s}\cdots a_1,a_{j_{s-1}}\cdots a_1\}\subseteq T$,
and $$a_{i_s}=(a_{i_s}\cdots a_1)\cdot (a_{j_{s-1}}\cdots a_1)\cdot a_{j_{s-1}+1}\cdot \ \cdots \ \cdot a_{i_s-1}\in \langle T\rangle.$$
If $i_{s}<j_s$, then $\{a_{j_{s-1}+1},\ldots,a_{i_s-1},a_{i_s+1},\ldots,a_{j_s}, a_{j_s}\cdots a_1,a_{j_{s-1}}\cdots a_1\}\subseteq T$, and so
$$a_{i_s}=(a_{j_s}\cdots a_1)\cdot (a_{j_{s-1}}\cdots a_1)\cdot a_{j_{s-1}+1}\cdot \ \cdots \ \cdot a_{i_s-1}\cdot a_{i_s+1}\cdot \ \cdots \ \cdot a_{j_s}\in \langle T\rangle.$$

Now, we have $a_{i_1},\ldots,a_{i_k}\in\langle T\rangle$. This means that $T$ is a generating subset of $\mz_2^n$, and since $|T|=n$,
it is a minimal generating subset.
It implies that for each sequence $(i_1,j_1,\ldots,i_k,j_k)$ satisfying inequation~\eqref{eq=2}, there exists a $Q_n$-isomorphic subgraph of $AQ_n$,
and since there are
$$f_1(n)=\sum\limits_{2\leq k \leq n-1}\
 \sum\limits_{2\leq j_1 < j_2<\cdots <j_k\leq n}{(j_k-j_{k-1})\cdots (j_2-j_1)j_1}$$
such sequences, we obtain $f_1(n)$ distinct $Q_n$-isomorphic subgraphs of $AQ_n$.

\medskip

In summary, there are
$$f(n)=1+\sum\limits_{2\leq i\leq n}{i}+f_1(n)=1+\sum\limits_{1\leq k \leq n-1}\
 \sum\limits_{2\leq j_1 < j_2<\cdots <j_k\leq n}{(j_k-j_{k-1})\cdots (j_2-j_1)j_1}$$
distinct subgraphs of $AQ_n$ that are isomorphic to $Q_n$.\hfill\qed

\begin{figure}
	\centering
	\includegraphics[width=0.5\linewidth]{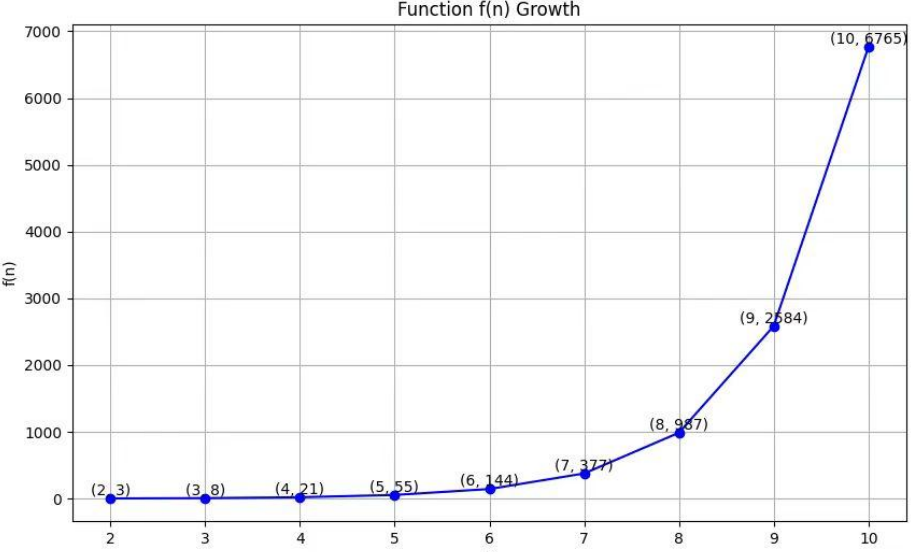}
	\caption{The function $f(n)$}
	\label{fig=3}
\end{figure}

Theorem~\ref{the=1} gives a lower bound of the number of $Q_n$-isomorphic subgraphs in $AQ_n$, which indicates that there are a lot of such subgraphs in $AQ_n$ (see Figure~\ref{fig=3}).
This bound is tight when $n=2$. It is easy to check that there are 3 different $Q_2$-isomorphic subgraphs in $AQ_2$, where
$3$ is equal to $f(2)$. When $n\geq3$, there remains a significant gap between this lower bound and the actual values. For example, we can find 45 different $Q_3$-isomorphic subgraphs in $AQ_3$, see Table~\ref{tab=1}, while $f(3)$ is only equal to 8.  However, in solving some practical problems related to $AQ_n$, we often do not require such a large number of $Q_n$-isomorphic subgraphs. This is  because many subgraphs exhibit only minor differences, with only a few edges distinguishing them from one another (see Rows 7 and 8 in Table~\ref{tab=1} for example). Identifying a few specific subgraphs is usually enough. The primary reason for employing the Cayley graphic method in this theorem to identify distinct $Q_n$-isomorphic subgraphs, lies in the fact that the Cayley graphic structure of $AQ_n$ allows us to easily find its subgraphs with common edges as few as possible. This proves to be highly convenient when addressing certain problems. In next subsection, we will present some of these subgraphs with least common edges, and provide applications in next two sections.

\begin{table}[ht]
{\tiny
\begin{center}
\begin{tabular}{l|l}
\hline
No.    & Edges whose induced subgraph is isomorphic to $Q_3$         \\
\hline
1  & $\{(000, 001), (001, 011), (011, 010), (010, 000), (100, 101),
    (101, 111), (111, 110), (110, 100), (000, 100), (001, 101),
    (010, 110), (011, 111)\}$ \\
\hline
2 & $\{(000, 001), (001, 011), (011, 010), (010, 000), (111, 101),
    (101, 100), (100, 110), (110, 111), (000, 111), (001, 101),
    (010, 110), (011, 100)\}$\\
\hline
3 &$\{(000, 010), (010, 011), (011, 001), (001, 000), (100, 101),
    (101, 111), (111, 110), (110, 100), (000, 100), (010, 101),
    (001, 110), (011, 111)\}$\\
\hline
4 & $\{(011, 001), (001, 000), (000, 010), (010, 011), (100, 101),
    (101, 111), (111, 110), (110, 100), (011, 100), (001, 101),
    (010, 110), (000, 111) \}$\\
\hline
5 & $\{(000, 001), (001, 011), (011, 010), (010, 000), (100, 110),
    (110, 111), (111, 101), (101, 100), (000, 100), (001, 110),
    (010, 101), (011, 111)\}$\\
\hline
6 & $\{(000, 011), (011, 010), (010, 001), (001, 000), (111, 100),
    (100, 110), (110, 101), (101, 111), (000, 111), (011, 100),
    (001, 101), (010, 110)\}$\\
\hline
7 & $\{(000, 011), (011, 001), (001, 010), (010, 000), (111, 100),
    {\bf (100, 110)},(110, 101), {\bf(101, 111)}, {\bf(000, 111)}, {\bf(011, 100)},
    (010, 101), (001, 110)\}$\\
\hline
8 & $\{(000, 011), (011, 001), (001, 010), (010, 000), (111, 100),
    {\bf(111, 110)}, (110, 101), {\bf(101, 100)}, {\bf(000, 100)}, {\bf(011, 111)},
    (010, 101), (001, 110)\}$\\
\hline
9 & $\{(011, 000), (000, 010), (010, 001), (001, 011), (111, 100),
    (100, 110), (110, 101), (101, 111), (011, 111), (000, 100),
    (001, 101), (010, 110)\}$\\
\hline
10 & $\{(000, 011), (011, 010), (010, 001), (001, 000), (111, 100),
    (100, 101), (101, 110), (110, 111), (000, 111), (011, 100),
    (001, 110), (010, 101)\}$\\
\hline
11 & $\{(000, 111), (111, 101), (101, 010), (010, 000), (100, 011),
    (011, 001), (001, 110), (110, 100), (000, 100), (111, 011),
    (010, 110), (101, 001) \}$\\
\hline
12 & $\{(000, 111), (111, 101), (101, 001), (001, 000), (100, 011),
    (011, 010), (010, 110), (110, 100), (000, 100), (111, 011),
    (001, 110), (101, 001)\}$\\
\hline
13 & $\{(000, 100), (100, 101), (101, 010), (010, 000), (111, 011),
    (011, 001), (001, 110), (110, 111), (000, 111), (100, 011),
    (010, 110), (101, 001)\}$\\
\hline
14 & $\{(011, 111), (111, 101), (101, 010), (010, 011), (100, 000),
    (000, 001), (001, 110), (110, 100), (011, 100), (111, 000),
    (010, 110), (101, 001)\}$\\
\hline
15 & $\{(000, 111), (111, 110), (110, 010), (010, 000), (100, 011),
    (011, 001), (001, 101), (101, 100), (000, 100), (111, 011),
    (010, 101), (110, 001)\}$\\
\hline
16 & $\{(000, 111), (111, 100), (100, 011), (011, 000), (010, 101),
    (101, 110), (110, 001), (001, 010), (000, 010), (111, 101),
    (011, 001), (100, 110) \}$\\
\hline
17 & $\{(000, 111), (111, 100), (100, 011), (011, 000), (001, 101),
    (101, 110), (110, 010), (010, 001), (000, 001), (111, 101),
    (011, 010), (100, 110)\}$\\
\hline
18 & $\{(000, 100), (100, 111), (111, 011), (011, 000), (010, 101),
    (101, 110), (110, 001), (001, 010), (000, 010), (100, 101),
    (011, 001), (111, 110)\}$\\
\hline
19 & $\{(011, 111), (111, 100), (100, 000), (000, 011), (010, 101),
    (101, 110), (110, 001), (001, 010), (011, 010), (111, 101),
    (000, 001), (100, 110)\}$\\
\hline
20 & $\{(000, 111), (111, 100), (100, 011), (011, 000), (010, 110),
    (110, 101), (101, 001), (001, 010), (000, 010), (111, 110),
    (011, 001), (100, 101)\}$\\
\hline
21 & $\{(000, 100), (100, 111), (111, 011), (011, 000), (001, 110),
    (110, 101), (101, 010), (010, 001), (000, 001), (100, 110),
    (011, 010), (111, 101)\}$\\
\hline
22 & $\{(000, 100), (100, 111), (111, 011), (011, 000), (010, 110),
    (110, 101), (101, 001), (001, 010), (000, 010), (100, 110),
    (011, 001), (111, 101)\}$\\
\hline
23 & $\{(000, 111), (111, 100), (100, 011), (011, 000), (001, 110),
    (110, 101), (101, 010), (010, 001), (000, 001), (111, 110),
    (011, 010), (100, 101)\}$\\
\hline
24 & $\{(011, 100), (100, 111), (111, 000), (000, 011), (001, 110),
    (110, 101), (101, 010), (010, 001), (011, 001), (100, 110),
    (000, 010), (111, 101)\}$\\
\hline
25 & $\{(000, 100), (100, 111), (111, 011), (011, 000), (001, 101),
    (101, 110), (110, 010), (010, 001), (000, 001), (100, 101),
    (011, 010), (111, 110) \}$\\
\hline
26 & $\{(000, 001), (001, 110), (110, 111), (111, 000), (010, 011),
    (011, 100), (100, 101), (101, 010), (000, 010), (001, 011),
    (111, 101), (110, 100)\}$\\
\hline
27 & $\{(000, 010), (010, 110), (110, 111), (111, 000), (001, 011),
    (011, 100), (100, 101), (101, 001), (000, 001), (010, 011),
    (111, 101), (110, 100)\}$\\
\hline
28 & $\{(000, 001), (001, 110), (110, 100), (100, 000), (010, 011),
    (011, 111), (111, 101), (101, 010), (000, 010), (001, 011),
    (100, 101), (110, 111)\}$\\
\hline
29 & $\{(011, 001), (001, 110), (110, 111), (111, 011), (010, 000),
    (000, 100), (100, 101), (101, 010), (011, 010), (001, 000),
    (111, 101), (110, 100)\}$\\
\hline
  30 &$\{
    (000, 001), (001, 101), (101, 111), (111, 000), (010, 011),
    (011, 100), (100, 110), (110, 010), (000, 010), (001, 011),
    (111, 110), (101, 100)\}$ \\
  \hline
  31 &$\{
    (000, 100), (100, 101), (101, 111), (111, 000), (001, 110),
    (110, 010), (010, 011), (011, 001), (000, 001), (100, 110),
    (111, 011), (101, 010) \}$\\
  \hline
  32 &$\{
    (000, 100), (100, 101), (101, 111), (111, 000), (010, 110),
    (110, 001), (001, 011), (011, 010), (000, 010), (100, 110),
    (111, 011), (101, 001) \}$\\
  \hline
  33 &$\{
    (000, 111), (111, 101), (101, 100), (100, 000), (010, 110),
    (110, 001), (001, 011), (011, 010), (000, 010), (111, 110),
    (100, 011), (101, 001)\}$ \\
  \hline
  34 &$\{
    (011, 100), (100, 101), (101, 111), (111, 011), (001, 110),
    (110, 010), (010, 000), (000, 001), (011, 001), (100, 110),
    (111, 000), (101, 010)\}$ \\
  \hline
  35 &$\{
    (000, 100), (100, 110), (110, 111), (111, 000), (001, 101),
    (101, 010), (010, 011), (011, 001), (000, 001), (100, 101),
    (111, 011), (110, 010)\}$ \\
  \hline
  36 &$\{
    (000, 100), (100, 110), (110, 001), (001, 000), (111, 011),
    (011, 010), (010, 101), (101, 111), (000, 111), (100, 011),
    (001, 101), (110, 010)\}$ \\
  \hline
  37 &$\{
    (000, 100), (100, 110), (110, 010), (010, 000), (111, 011),
    (011, 001), (001, 101), (101, 111), (000, 111), (100, 011),
    (010, 101), (110, 001)\}$ \\
  \hline
  38 &$\{
    (000, 111), (111, 110), (110, 001), (001, 000), (100, 011),
    (011, 010), (010, 101), (101, 100), (000, 100), (111, 011),
    (001, 101), (110, 010)\}$ \\
  \hline
  39 &$\{
    (011, 100), (100, 110), (110, 001), (001, 011), (111, 000),
    (000, 010), (010, 101), (101, 111), (011, 111), (100, 000),
    (001, 101), (110, 010)\}$ \\
  \hline
  40 &$\{
    (000, 100), (100, 101), (101, 001), (001, 000), (111, 011),
    (011, 010), (010, 110), (110, 111), (000, 111), (100, 011),
    (001, 110), (101, 010)\}$ \\
  \hline
  41 &$\{
    (000, 100), (100, 011), (011, 010), (010, 000), (111, 101),
    (101, 001), (001, 110), (110, 111), (000, 111), (100, 101),
    (010, 110), (011, 001)\}$ \\
  \hline
  42 &$\{
    (000, 100), (100, 011), (011, 001), (001, 000), (111, 101),
    (101, 010), (010, 110), (110, 111), (000, 111), (100, 101),
    (001, 110), (011, 010)\}$ \\
  \hline
  43 &$\{
    (000, 111), (111, 011), (011, 010), (010, 000), (100, 101),
    (101, 001), (001, 110), (110, 100), (000, 100), (111, 101),
    (010, 110), (011, 001)\}$ \\
  \hline
  44 &$\{
    (011, 100), (100, 000), (000, 010), (010, 011), (111, 101),
    (101, 001), (001, 110), (110, 111), (011, 111), (100, 101),
    (010, 110), (000, 001)\}$ \\
  \hline
  45 & $\{
    (000, 100), (100, 011), (011, 010), (010, 000), (111, 110),
    (110, 001), (001, 101), (101, 111), (000, 111), (100, 110),
    (010, 101), (011, 001)\}$ \\
\hline
\end{tabular}
\end{center}
}
\vskip -0.5cm
\caption{ The edge sets of $Q_3$-isomorphic subgraphs in $AQ_3$}\label{tab=1}
\end{table}

\subsection{Two $Q_n$-isomorphic subgraphs of $AQ_n$ with least number of common edges}

Let $G_1$ and $G_2$ be two different $Q_n$-isomorphic subgraphs of $AQ_n$.
Then $E(G_1)\cup E(G_2)\subseteq E(AQ_n)$.
Noting that $AQ_n$ is a $(2n-1)$-regular graph
while both $G_1$ and $G_2$ are $n$-regular graphs, we have that $|E(AQ_n)|=(2n-1)2^{n-1}$ and $|E(G_1)|=|E(G_2)|=n2^{n-1}$. This forces that the graphs $G_1$ and $G_2$ must have common edges, i.e., $E(G_1)\cap E(G_2)\neq \emptyset$. Moreover, we have
\begin{equation}\label{eq=3.2}
|E(G_1)\cap E(G_2)|=|E(G_1)|+|E(G_2)|-|E(G_1)\cup E(G_2)|\geq |E(G_1)|+|E(G_2)|-|E(AQ_n)|\geq 2^{n-1}.
\end{equation}
In particular, if $ |E(G_1)\cap E(G_2)|=2^{n-1}$, then Eq.~\eqref{eq=3.2} forces that $E(G_1)\cup E(G_2)=E(AQ_n)$.

\begin{lem}\label{lem=3.2-1}
Let $G_1$ and $G_2$ be two different $Q_n$-isomorphic subgraphs of $AQ_n$. If
$|E(G_1)\cap E(G_2)|=2^{n-1}$, then $E(G_1)\cap E(G_2)$ is a perfect matching of $AQ_n$.
\end{lem}

\f {\bf Proof.} Consider the edge-induced subgraph of $E(G_1)\cap E(G_2)$ in $AQ_n$, denoted as $H$. Then $|E(H)|=2^{n-1}$.
Since $G_1$ and $G_2$ are two spanning subgraphs of $AQ_n$ that are isomorphic to $Q_n$, each vertex $v\in V(AQ_n)$ is incident to $n$ edges in both $G_1$ and $G_2$,
and since $d_{AQ_n}(v)=2n-1$,
at least one edge incident to $v$ belongs to $E(G_1)\cap E(G_2)$. It forces that $v\in V(H)$ and $d_{H}(v)\geq1$, i.e., $H$ is a spanning subgraph of $AQ_n$ of minimum degree at least 1.
Since $|V(H)|=2^n$ and $|E(H)|=2^{n-1}$, we have that $H$ is a 1-regular graph and so  $E(H)=E(G_1)\cap E(G_2)$ is a perfect matching of $AQ_n$.\hfill\qed


The subsequent construction will establish the existence of two $Q_n$-isomorphic subgraphs of $AQ_n$ whose edge-sets have intersection a perfect matching. We also use $\mathcal{CG}_n=\Cay(\mz_2^n,S)$
to represent the augmented cube $AQ_n$.
Let
$$\begin{array}{ll}
E_i&=\{(g,a_ig)~|~g\in\mz_2^n\},~1\leq i\leq n,\\
E_{\leq j}&=\{(g,(a_ja_{j-1}\cdots a_1)g)~|~g\in\mz_2^n\},~2\leq j\leq n.
  \end{array}$$
In the language of binary strings (see the isomorphism $\phi$ given in the proof of Lemma~\ref{lem=1}), these subsets above can be also written as follows
$$\begin{array}{ll}
E_i&=\{(x_n\cdots x_1,  x_n \cdots x_{i+1}(x_i+1)x_{i-1}\cdots x_1)\  | \ x_n\cdots x_1\in V(AQ_n)\},~1\leq i\leq n,\\
E_{\leq j}&=\{(x_n\cdots x_1,x_n \cdots x_{j+1}(x_j+1)\cdots (x_1+1))~|~x_n\cdots x_1\in V(AQ_n)\},~2\leq j\leq n.
\end{array}$$
By Definition~\ref{defi=1}, $E_i$ is the set of
all hypercube edges of dimension $i$ in $AQ_n$,
while $E_{\leq j}$ is the set of all augmented edges of dimension $j$. In particular, each of these $E_i$s and $E_{\leq j}$s is a perfect matching of $AQ_n$, and
$$E(AQ_n)=(\cup_{i=1}^n E_i) \cup (\cup_{j=2}^n E_{\leq j}).$$

Now, we partition the edge set of $E(AQ_n)$ into $2n-1$ perfect matchings. Considering the edge-induced subgraphs in $AQ_n$ induced by

\parbox{8cm}{
\begin{eqnarray*}
&& E_1\cup E_2\cup \cdots \cup E_n,\\
&& E_1\cup E_{\leq 2}\cup \cdots \cup  E_{\leq n}.
\end{eqnarray*}}\hfill
\parbox{1cm}{
\begin{eqnarray}
\label{eq=3}\\ \label{eq=4}
\end{eqnarray}}

\f Clearly, they are Cayley graph $\Cay(\mz_2^n,\{a_1,a_2,\ldots, a_n\})$ and $\Cay(\mz_2^n,\{a_1,a_2a_1,\ldots,a_na_{n-1}\cdots a_1\})$, respectively. By Proposition~\ref{prop=2.3-1}, both of them are isomorphic to $Q_n$. We denoted these two subgraphs as $Q_n^1$ and $Q_n^2$, respectively, which will be used frequently in later section.
Moreover, $E(Q_n^1)\cap E(Q_n^2)=E_1$.

\begin{figure}
	\centering
	\includegraphics[width=0.4\linewidth]{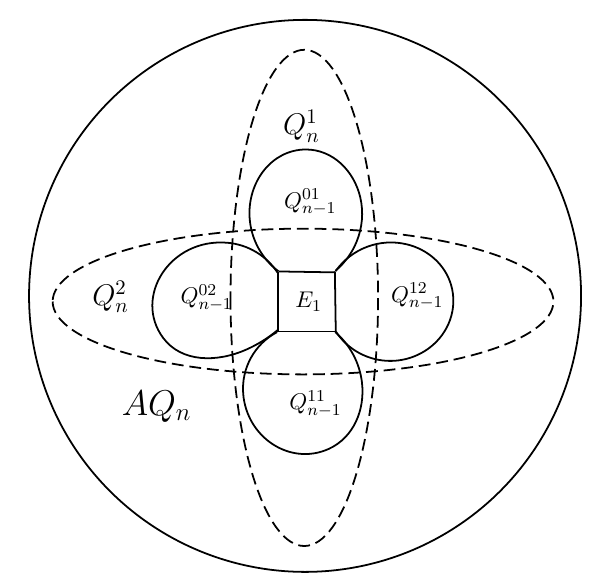}
	\caption{A decomposition of the edge set of $AQ_n$}
	\label{fig=2}
\end{figure}

It is well-known that in a hypercube $Q_n$, after deleting all $k$-dimensional edges, the remaining subgraph is disconnected and consists of two disjoint $(n-1)$-dimensional hypercubes~\cite{H1}. It is referred as a {\em partition} of $Q_n$ along dimension $k$, and often denoted as $Q_n=Q_{n-1}^0\odot Q_{n-1}^1$.
Hence, in $Q_n^1$ and $Q_n^2$, when we remove the edges belonging to $E_1$, four $Q_{n-1}$-isomorphic subgraphs are obtained. Since $E(Q_n^1)\cap E(Q_n^2)=E_1$, these four $Q_{n-1}$-isomorphic subgraphs are edge-disjoint. An illustration is presented in Figure~\ref{fig=2}.

In summary, we have the following lemma.

\begin{lem}\label{lem=2}
In $AQ_n$ with $n\geq2$, $Q_{n}^1$ and $Q_n^2$ are two $Q_n$-isomorphic subgraphs in
$AQ_n$ with $E(Q_n^1)\cap E(Q_n^2)=E_1$. Furthermore, there exists  four edge-disjoint $Q_{n-1}$-isomorphic subgraphs,
denoted as
$$Q^1_n=Q^{01}_{n-1}\odot Q^{11}_{n-1},\ Q^2_n=Q^{02}_{n-1}\odot Q^{12}_{n-1},$$
whose vertex sets are
$$\begin{array}{ll}
&V(Q^{01}_{n-1})=\{x_{n}\cdots x_20 \ | \ x_i\in \{0,1\}, 2\leq i\leq n\};\\
&V(Q^{11}_{n-1})=\{x_{n}\cdots x_21 \ | \ x_i\in \{0,1\}, 2\leq i\leq n\};\\
&V(Q^{02}_{n-1})=\{x_nx_{n-1}\cdots x_2x_2\ | \ x_i\in \{0,1\}, 2\leq i\leq n\};\\
&V(Q^{12}_{n-1})=\{x_nx_{n-1}\cdots x_2(1+x_2)\ | \ x_i\in \{0,1\}, 2\leq i\leq n\}.
\end{array}$$
\end{lem}

\f {\bf Proof.} We only need to explain the conclusion of $V(Q^{02}_{n-1})$
and $V(Q^{12}_{n-1})$. In $Q_n^2$, i.e., the Cayley graph $\Cay(\mz_2^n,\{a_1,a_2a_1,a_3a_2a_1,\ldots,a_na_{n-1}\cdots a_1\})$, after removing the edges in $E_1=\{(g,a_1g)~|~g\in\mz_2^n\}$, the remaining subgraph is disconnected, which is the Cayley graph
$\Cay(\mz_2^n,\{a_2a_1,a_3a_2a_1,\ldots,a_na_{n-1}\cdots a_1\})$. Since $\langle a_2a_1,a_3a_2a_1,\ldots,a_na_{n-1}\cdots a_1\rangle\cong\mz_2^{n-1}$, we have $|\mz_2^n|/|\langle a_2a_1,a_3a_2a_1,\ldots,a_na_{n-1}\cdots a_1\rangle|=2$, and so there are two components is $Q_n^2-E_1$ whose vertex sets are
$$\begin{array}{ll}
&\langle a_2a_1,a_3a_2a_1,\ldots,a_na_{n-1}\cdots a_1\rangle\\
=&\{(a_2a_1)^{y_2}(a_3a_2a_1)^{y_3}\cdots
(a_na_{n-1}\cdots a_1)^{y_n}~|~y_i\in\{0,1\},2\leq i\leq n\}\\
=&\{a_n^{y_n}a_{n-1}^{y_{n-1}+y_n}\cdots
a_{2}^{y_2+\cdots+y_{n-1}+y_n}a_{1}^{y_2+\cdots+y_{n-1}+y_n}~|~y_i\in\{0,1\},2\leq i\leq n\};\\
&\langle a_2a_1,a_3a_2a_1,\ldots,a_na_{n-1}\cdots a_1\rangle a_1\\
=&\{(a_2a_1)^{y_2}(a_3a_2a_1)^{y_3}\cdots
(a_na_{n-1}\cdots a_1)^{y_n}a_1~|~y_i\in\{0,1\},2\leq i\leq n\}\\
=&\{a_n^{y_n}a_{n-1}^{y_{n-1}+y_n}\cdots
a_{2}^{y_2+\cdots+y_{n-1}+y_n}a_{1}^{1+y_2+\cdots+y_{n-1}+y_n}~|~y_i\in\{0,1\},2\leq i\leq n\},
\end{array}$$
respectively. We denoted these two components as $Q_{n-1}^{02}$ and $Q_{n-1}^{12}$, respectively.
In the language of binary strings (see the isomorphism $\phi$ given in the proof of Lemma~\ref{lem=1}), the vertex sets of $Q_{n-1}^{02}$ and $Q_{n-1}^{12}$ can be written as follows:
$$\begin{array}{ll}
V(Q^{02}_{n-1})&=\{y_n(y_{n-1}+y_n)\cdots (y_2+\cdots +y_n)(y_2+\cdots +y_n)\ | \ y_i\in \{0,1\}, 2\leq i\leq n\}\\
&=\{x_nx_{n-1}\cdots x_2x_2\ | \ x_i\in \{0,1\}, 2\leq i\leq n\};\\
V(Q^{12}_{n-1})&=\{y_n(y_{n-1}+y_n)\cdots (y_2+\cdots +y_n)(1+y_2+\cdots +y_n)\ | \ y_i\in \{0,1\}, 2\leq i\leq n\}\\
&=\{x_nx_{n-1}\cdots x_2(1+x_2)\ | \ x_i\in \{0,1\}, 2\leq i\leq n\}.
\end{array}$$
The proof is complete.\hfill\qed

Let
\begin{equation}\label{eq=3.1}
\varphi_1: x_n\cdots x_20 \mapsto x_{n}\cdots x_21, \ \forall \ x_{n}\cdots x_20 \in V(Q_{n-1}^{01})
\end{equation}
and
\begin{equation}\label{eq=3.2}
\varphi_2:x_nx_{n-1}\cdots x_2x_2 \mapsto x_nx_{n-1}\cdots x_2(1+x_2), \ \forall \ x_nx_{n-1}\cdots x_2x_2 \in V(Q_{n-1}^{02}).
\end{equation}
It is easy to find that $\varphi_i$ is an isomorphism from $Q^{0i}_{n-1}$ to $Q^{1i}_{n-1}$ for $i=1$ or 2.
Based on Lemmas~\ref{lem=3.2-1} and~\ref{lem=2},
we have the following theorem.

\begin{theorem}
Let $G_1$ and $G_2$ be two of the $Q_n$-isomorphic subgraphs of $AQ_n$ such that they have the minimum number of common edges. Then $E(G_1)\cap E(G_2)$ is a perfect matching of $AQ_n$.
\end{theorem}

We note that the pair of $Q_n$-isomorphic subgraphs of $AQ_n$ with the minimum number of common edges is not unique. Based on the Cayley property of $AQ_n$, the following part provides more such subgraphs.

Let $k$ be an integer with $1\leq k\leq n-1$. For any $k$ integers $j_1,\ldots,j_k$
with $2\leq j_1<j_2<\cdots <j_k\leq n$. It can be checked easily that the following two subsets are both minimal generalized subsets of $\mz_2^n$:
$$\begin{array}{ll}
S_3=\{&a_1,a_2,\ldots, a_{j_1-1},a_{j_1}\cdots a_{1},a_{j_1+1},\ldots,a_{j_2-1},a_{j_2}\cdots a_{1},a_{j_2+1},\ldots, a_{j_k-1},a_{j_k}\cdots a_{1},a_{j_k+1},\ldots,a_n\};\\
S_4=\{&a_1,a_2a_1,\ldots, a_{j_1}\cdots a_{1},a_{j_1},a_{j_1+1}\cdots a_{1},\ldots,a_{j_2-1}\cdots a_{1},a_{j_2},a_{j_2+1}\cdots a_{1},\ldots, a_{j_k-1}\cdots a_1, a_{j_k},\\
&a_{j_k+1}\cdots a_1,\ldots,a_n\cdots a_1\}.
\end{array}$$
These are indeed the subsets through exchanging $a_{j_1},a_{j_2},\ldots,a_{j_k}$ in
$\{a_1,a_2,\ldots, a_n\}$
with the elements $a_{j_1}\cdots a_1, a_{j_2}\cdots a_2,\ldots,a_{j_k}\cdots a_1$ in $\{a_1,a_2a_1,\ldots, a_n\cdots a_1\}$.
Therefore, by Proposition~\ref{prop=2.3-1} the Cayley graphs of $\mz_2^n$ with respect to $S_3$ and $S_4$ are both isomorphic to $Q_n$, and so we have a more general conclusion than Lemma~\ref{lem=2}.

\begin{theorem}\label{the=3.2-1}
Let $k$ be an integer with $1\leq k\leq n-1$. For any $k$ integers $j_1,\ldots,j_k$
with $2\leq j_1<j_2<\cdots <j_k\leq n$, the subgraphs induced by the two subsets of edges
{\small $$\begin{array}{ll}
&E_1\cup E_2\cup \cdots \cup E_{j_1-1}\cup E_{\leq j_1}\cup E_{j_1+1}\cup \cdots\cup E_{j_2-1}\cup E_{\leq j_2}\cup E_{j_2+1}\cup \cdots\cup E_{j_k-1}\cup E_{\leq j_k}\cup E_{j_k+1}\cup\cdots \cup E_{n};\\
&E_1\cup E_{\leq 2}\cup \cdots \cup E_{\leq j_1-1}\cup E_{j_1}\cup E_{\leq j_1+1}\cup \cdots\cup E_{\leq j_2-1}\cup E_{j_2}\cup E_{\leq j_2+1}\cup \cdots\cup E_{\leq j_k-1}\cup E_{j_k}\cup E_{\leq j_k+1}\cup\cdots \cup E_{\leq n},
\end{array}$$}
in $AQ_n$ are both isomorphic to $Q_n$.
\end{theorem}

\begin{figure}
	\centering
	\includegraphics[width=1\linewidth]{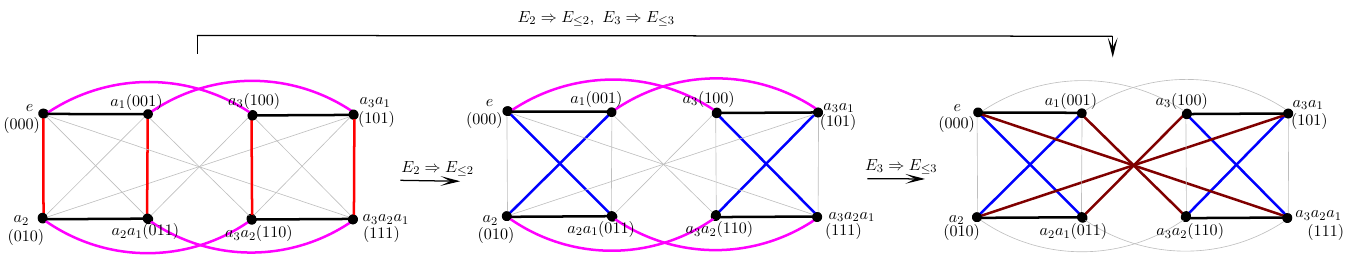}
	\caption{Construction of $Q_3$ in $AQ_3$ via PM-Reciprocity Construction Method}
	\label{fig=3.5}
\end{figure}

In Theorem~\ref{the=3.2-1}, the construction of a pair of $Q_n$-isomorphic subgraphs in $AQ_n$ is realized through a reciprocal exchange of perfect matchings between $Q_n^1$ and $Q_n^2$. An example is given in Figure~\ref{fig=3.5}. The selection of these mutually exchanged perfect matchings depends on the choice of a minimal generating set for $\mz_2^n$. We term this technique the {\em Perfect Matching Reciprocal Construction Method} or {\em PM-Reciprocity Construction Method}. Notably, the edge sets of the pair of subgraphs constructed above share $E_1$ as their common intersection. Furthermore, by applying the same method,
we can construct two $Q_n$-isomorphic subgraphs
in $AQ_n$, whose edge-sets have intersection
$E_j$ or $E_{\leq j}$ for any $2\leq j\leq n$.

\begin{theorem}\label{the=3.2-2}
For any integer $j$
with $2\leq j\leq n$,
there exist two $Q_n$-isomorphic subgraphs in $AQ_n$ such that the intersection of their edge-sets is equal to $E_j$ or $E_{\leq j}$.
\end{theorem}

\f {\bf Proof.} We shall proceed to prove this theorem from the perspective of Cayley graphs, i.e., by letting $AQ_n=\mathcal{CG}_n=\Cay(\mz_2^n, S)$.
For $2\leq j\leq n$, considering the following two edge-induced subgraphs

\parbox{8cm}{
\begin{eqnarray*}
&& G_1:=AQ_n[E_2\cup \cdots \cup E_n\cup E_{\leq j}],\\
&& G_2:=AQ_n[E_1\cup E_{\leq 2}\cup \cdots\cup E_{\leq j-1}\cup E_j\cup E_{\leq j+1}\cup \cdots\cup E_{\leq n}].
\end{eqnarray*}}\hfill
\parbox{1cm}{
\begin{eqnarray}
\label{eq=5}\\ \label{eq=6}
\end{eqnarray}}

\f (For the case where $j=2$, $E_{\leq j-1}$ is regarded as $E_1$.) Then
$$
\begin{array}{ll}
& G_1=\Cay(\mz_2^n,\{a_2,\ldots,a_n,a_{j}\cdots a_1\});\\
& G_2=\Cay(\mz_2^n,\{a_1,a_2a_1,\ldots,a_{j-1}\cdots a_1, a_j, a_{j+1}\cdots a_1,\ldots,a_{n}\cdots a_1\}).
\end{array}
$$
Clearly, both $\{ a_2,\ldots,a_j,a_{j}\cdots a_1,a_{j+1},\ldots, a_n\}$
and $\{a_1,a_2a_1,\ldots,a_{j-1}\cdots a_1, a_j, a_{j+1}\cdots a_1,\ldots,a_{n}\cdots a_1\}$ are minimal generating subsets of $\mz_2^n$, and so by Proposition~\ref{prop=2.3-1} both $G_1$ and $G_2$ are isomorphic to $Q_n$, and $E(G_1)\cap E(G_2)=E_j$. Indeed, the subgraph $G_1$ (or $G_2$, resp.) is derived from $Q_n^1$ (or $Q_n^2$, resp.) by replacing $E_1$ (or $E_{\leq j}$, resp.) with $E_{\leq j}$ (or $E_{j}$, resp.).

Recall that $E(Q_n^2)=E_1\cup E_{\leq 2}\cup \cdots \cup E_{\leq n}$. Then $E(Q_n^2)\cap E(G_1)=E_{\leq j}$. In other words, $Q_n^2$
and $G_1$ are two $Q_n$-isomorphic subgraphs of $AQ_n$ whose edge-sets have intersection $E_{\leq j}$.
\hfill\qed

The above two theorems, together with the following lemma, will be used frequently in the following proof. Since
$\{a_2,\ldots, a_{j_1},a_{j_1}\cdots a_1,a_{j_1+1},\ldots,a_{j_2-1},a_{j_2}\cdots a_1,a_{j_2+1},\ldots, a_{j_n}\}$
is a minimal generating subset of $\mz_2^n$, where $j_1$ and $j_2$ are two integers satisfying $2\leq j_1<j_2\leq n$, we have the following lemma.

\begin{lem}\label{lem=3.2-2}
For any two integers $j_1$ and $j_2$ with $2\leq j_1<j_2\leq n$, the edge-induced subgraph of $AQ_n$ induced by
$E_2\cup \cdots \cup E_{j_1}\cup E_{\leq j_1}\cup E_{j_1+1}\cup \cdots \cup E_{j_2-1}\cup E_{\leq j_2}\cup E_{j_2+1}\cup \cdots \cup E_n$
is isomorphic to $Q_n$.
\end{lem}

\section{EDHCs in an augmented cube}

In this section, we aim to construct the edge-disjoint Hamiltonian cycles (EDHCs) in an augmented cube $AQ_n$. The following proposition related to hypercubes comes from~\cite{BB}.

\begin{prop}\label{prop=4.1}
A hypercube $Q_n$ has $\lfloor \frac{n}{2} \rfloor$ EDHCs when $n\geq 2$.
\end{prop}

\begin{theorem}\label{theo=4.1}
There are $m$ EDHCs in $AQ_n$ with $n\geq 3$, where
$m=n-1$ when $n$ is odd and $m=n-2$ when $n$ is even.
\end{theorem}

\f {\bf Proof.} There are two EDHCs in $AQ_3$ and $AQ_4$ by~\cite[Theorem 6]{Hung1}.
Therefore, the theorem holds when $n=3$ or $4$. Assume $n\geq5$ in the following proof. We partition $AQ_n$ to four edge-disjoint subgraphs. By Lemma~\ref{lem=2}, there exist four
edge-disjoint $Q_{n-1}$-isomorphic subgraphs $Q_{n-1}^{01}$, $Q_{n-1}^{11}$, $Q_{n-1}^{02}$
and $Q_{n-1}^{12}$. The $Q_n^1=Q^{01}_{n-1}\odot Q^{11}_{n-1}$
and $Q^2_n=Q^{02}_{n-1}\odot Q^{12}_{n-1}$
are two spanning subgraphs of $AQ_n$, whose edge sets have intersection $E_1$, i.e. the set of $1$-dimensional hypercube edges (see Figure~\ref{fig=2} for an illustration).
The main idea of this proof is to construct $\lfloor \frac{n-1}{2}\rfloor$ EDHCs separately within these two spanning subgraphs.

Firstly, in $Q_{n}^{0}=Q^{01}_{n-1}\odot Q^{11}_{n-1}$, since $Q_{n-1}^{01}\cong Q_{n-1}$, it has $\lfloor \frac{n-1}{2}\rfloor$ EDHCs by Proposition~\ref{prop=4.1}, denoted by $C_1,\ldots, C_{\left\lfloor\frac{n-1}{2}\right\rfloor}$. Note that each of these cycles has length $2^{n-1}$. Since $n\geq 5$, we have $2^{n-1}-4(\left\lfloor\frac{n-1}{2}\right\rfloor -1) > 0$. It yields that there are $\lfloor \frac{n-1}{2}\rfloor$ edges $(u_1, v_1),\ldots , (u_{\left\lfloor\frac{n-1}{2}\right\rfloor}, v_{\left\lfloor\frac{n-1}{2}\right\rfloor})$, which lie on $C_1,\ldots, C_{\left\lfloor\frac{n-1}{2}\right\rfloor}$ respectively, and any two of them have no common end-vertices.
Since $\varphi_1$ is an isomorphism from $Q^{01}_{n-1}$ to $Q^{11}_{n-1}$ (see Eq.~\eqref{eq=3.1}), we have that $\varphi_1(C_1),\ldots , \varphi_1(C_{\left\lfloor\frac{n-1}{2}\right\rfloor})$ are EDHCs of $Q^{11}_{n-1}$. In particular, $\varphi_1(u_i,v_i)$ is an edge of $\varphi_1(C_i)$
for $1\leq i\leq \lfloor \frac{n-1}{2} \rfloor$, and any two of these edges have no common end-vertices.
In view of Eq.~\eqref{eq=3.1} and Definition~\ref{defi=1}~(1), for any $1\leq i\leq \left\lfloor\frac{n-1}{2}\right\rfloor$, both $(u_i, \varphi_1 (u_i))$ and $(v_i, \varphi_1 (v_i))$ are $1$-dimensional hypercube edges.
Now, we let
$$\mathbf{C}_i=C_i-(u_i, v_i)+(\varphi_1 (u_i), u_i)+\varphi_1(C_i)-(\varphi_1 (u_i), \varphi_1 (v_i))+(\varphi_1 (v_i), v_i)).$$
Then $\mathbf{C}_i$ is a Hamiltonian cycle of $Q^1_n=Q^{01}_{n-1}\oplus Q^{11}_{n-1}$, and so is in $AQ_n$.
For any $1\leq i\neq j\leq \left\lfloor\frac{n-1}{2}\right\rfloor$,
since $\{u_i,v_i\}\cap \{u_j,v_j\}=\{\varphi_1(u_i),\varphi_1(v_i)\}\cap \{\varphi_1(u_j),\varphi_1(v_j)\}=\emptyset$,
we have
$$\{(u_i, \varphi_1 (u_i)), \ (v_i, \varphi_1 (v_i))\}\cap \{(u_j, \varphi_1 (u_j)), \ (v_j, \varphi_1 (v_j))\}=\emptyset,$$
and since
$C_i,C_j, \varphi_1(C_i),\varphi_1(C_j)$ are edge-disjoint, we have that $\mathbf{C}_1,
\mathbf{C}_2,\ldots,\mathbf{C}_{\lfloor \frac{n-1}{2}\rfloor}$ are edge-disjoint.
Consequently, in the subgraph $Q_n^1=Q_{n-1}^{01}\odot Q_{n-1}^{11}$ we construct
$\lfloor \frac{n-1}{2}\rfloor$ EDHCs of $AQ_n$.

\medskip

Next, we turn to $Q^2_{n-1}=Q^{02}_{n-1}\odot Q^{12}_{n-1}$, and construct other $\lfloor \frac{n-1}{2}\rfloor$ EDHCs in the same way.
Again by Proposition~\ref{prop=4.1}, $Q^{02}_{n-1}$ has $\left\lfloor\frac{n-1}{2}\right\rfloor$ EDHCs, denoted by $D_1,\ldots , D_{\left\lfloor\frac{n-1}{2}\right\rfloor}$.
Note that during the construction of Hamiltonian cycles
in $Q_{n}^{1}$, we have already utilized $2\cdot \left\lfloor\frac{n-1}{2}\right\rfloor$ vertices $u_1,v_1,\ldots,u_{\lfloor\frac{n-1}{2}\rfloor},v_{\lfloor\frac{n-1}{2}\rfloor}$. Moving forward, our objective is to select distinct vertices in $D_1,\ldots , D_{\left\lfloor\frac{n-1}{2}\right\rfloor}$, to ensure that the cycles finally constructed are edge-disjoint with $\mathbf{C}_1,
\mathbf{C}_2,\ldots,\mathbf{C}_{\lfloor \frac{n-1}{2}\rfloor}$.

Since $2^{n-1}-4\left\lfloor\frac{n-1}{2}\right\rfloor -4(\left\lfloor\frac{n-1}{2}\right\rfloor-1) > 0$ when $n\geq5$, there are edges $(x_1, y_1),\ldots ,(x_{\left\lfloor\frac{n-1}{2}\right\rfloor},\  y_{\left\lfloor\frac{n-1}{2}\right\rfloor})$ on $D_1,\ldots ,D_{\lfloor\frac{n-1}{2}\rfloor}$ respectively, such that for any $1\leq i, j, k\leq \lfloor\frac{n-1}{2}\rfloor$ with $i\neq j$,
\begin{equation}\label{eq=4.1}
\{x_i, y_i\}\cap \{x_j, y_j\}= \{x_i, y_i\}\cap \{u_k, v_k\}=\emptyset.
\end{equation}
Since $\varphi_2$ is an isomorphism from $Q_{n-1}^{02}$ to $Q_{n-1}^{12}$ (see Eq.~\eqref{eq=3.2}), we have that
$\varphi_2(D_1),\ldots , \varphi_2(D_{\lfloor\frac{n-1}{2}\rfloor})$ are edge-disjoint Hamiltonian cycles of $Q^{12}_{n-1}$. In addition, $\varphi_2(x_1, y_1),\ldots ,\varphi_2 (x_{\lfloor\frac{n-1}{2}\rfloor}, y_{\lfloor\frac{n-1}{2}\rfloor})$  are edges of $\varphi_2(D_1),\ldots , \varphi_2 (D_{\lfloor\frac{n-1}{2}\rfloor})$, and for any $1\leq i\neq j,\ k\leq \lfloor\frac{n-1}{2}\rfloor$,
$$\{\varphi_2(x_i),\  \varphi_2(y_i)\}\cap \{\varphi_2(x_j),\ \varphi_2( y_j)\}= \{\varphi_2 (x_i), \ \varphi_2(y_i)\}\cap \{\varphi_2(u_k),\ \varphi_2( v_k)\}=\emptyset .$$
Let
$$\mathbf{D}_i=D_i-(x_i, y_i)+(x_i, \varphi_2 (x_i))+\varphi_2(D_i)-(\varphi_2(x_i), \varphi_2 (y_i))+(\varphi_2(y_i), y_i)),~1\leq i\leq \lfloor\frac{n-1}{2}\rfloor,$$
Then $\mathbf{D}_1,\ldots , \mathbf{D}_{\lfloor\frac{n-1}{2}\rfloor}$  are EDHCs of $Q^2_n$ (also of $AQ_n$).

Recall that $E(Q^1_n)\cap E(Q^2_n)=E_1$.
For any $1\leq i,k\leq \lfloor\frac{n-1}{2}\rfloor$,
the cycle $\mathbf{C}_i$ has only two $1$-dimensional hypercube edges that are $(u_i,\varphi_1(u_i))$ and $(v_i,\varphi_1(v_i))$,
while $\mathbf{D}_i$ also has only two $1$-dimensional hypercube edges that are $(x_i,\varphi_2(x_i))$ and $(y_i,\varphi_2(y_i))$.
In view of Eq.~\eqref{eq=4.1}, these $1$-dimensional edges are disjoint.
It implies that $\mathbf{C}_1,\ldots , \mathbf{C}_{\left\lfloor\frac{n-1}{2}\right\rfloor}, \mathbf{D}_1,\ldots , \mathbf{D}_{\lfloor\frac{n-1}{2}\rfloor}$ are edge-disjoint.

Finally, we have constructed $2\cdot \lfloor\frac{n-1}{2}\rfloor$ EDHCs of $AQ_n$.
When $n$ is odd, $2\cdot \left\lfloor\frac{n-1}{2}\right\rfloor=2\cdot \frac{n-1}{2}=n-1$; and when $n$ is even, $2\cdot {\left\lfloor\frac{n-1}{2}\right\rfloor}=2\cdot (\frac{n}{2}-1)=n-2$. The proof is complete.
\hfill\qed

Noting that $\lfloor\frac{2n-1}{2}\rfloor=n-1$, Theorem~\ref{theo=4.1} implies that there are $\lfloor\frac{2n-1}{2}\rfloor$ EDHCs in $AQ_n$ when $n$ is odd.
Since $AQ_n$ is a $(2n-1)$-regular graph, this confirms Conjecture~\ref{con=1} when $n$ is odd.

\begin{cor}
  The $AQ_n$ admits a Hamiltonian decomposition when $n\geq3$ is odd.
\end{cor}

\section{Edge-fault-tolerant Hamiltonicity of $AQ_n$}

In this section, we consider the fault-tolerant Hamiltonicity of the augmented cube.
The following proposition comes from~\cite[Theorem 4.1]{XZZY}.

\begin{prop}\label{prop=5.3-1}
In $AQ_n$ with $n\geq3$, let $A\subset V(AQ_n)\cup E(AQ_n)$ with $|A|\leq 2n-4$, and let $f$ be the number of vertices belonging to $A$. For any distinct vertices $u$ and $v$ in $AQ_n-A$,
there exists a $(u,v)$-path of each length from
$\max\{d_{AQ_n}(u,v)+2,4\}$ to $2^n-f-1$.
\end{prop}

\begin{figure}
	\centering	\includegraphics[width=0.2\linewidth]{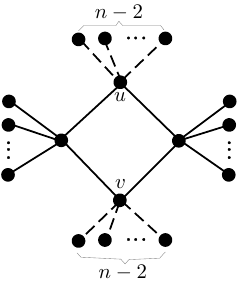}
	\caption{An example of $Q_n$ with $2n-4$ faulty edges (In this illustration, fault-free edges are depicted as solid lines, and faulty edges as dashed lines; all following figures adhere to this convention.)}
	\label{fig=4}
\end{figure}

In $AQ_n$ or $Q_n$, we assume that all vertices are fault-free,
and use $F$ to denote the set of faulty edges in this section. Chan and Lee~\cite{CH} proved that in a faulty hypercube $Q_n$ with $n\geq 3$,  if $|F|\leq 2n-5$ and
each vertex is incident to at least two fault-free edges, then $Q_n-F$ has a Hamiltonian cycle. This result is optimal in term of the number of faulty edges in $Q_n$, since there exists a hypercube with $2n-4$ faulty edges in which no fault-free Hamiltonian cycle exists.
An example of $Q_n$ with $n\geq3$ and $2n-4$ faulty edges is illustrated in Figure~\ref{fig=4},
in which $u$ and $v$ are two non-adjacent vertices belonging to a fault-free $4$-cycle $C$ and all edges incident to $u$ and $v$ in $Q_n$, expect those in $C$, are faulty. Subsequent research has significantly extended these findings. Later, Liu and Wang~\cite{LW} proved that under the conditional fault model, Chan and Lee's conclusion remains valid for up to $3n-8$ faulty edges, excluding the case depicted in Figure~\ref{fig=4}. Furthermore, we have a stronger result as follows, coming from~\cite[Theorem~1]{CH}.

\begin{prop}\label{prop=5.2}
In a $($faulty$)$ hypercube $Q_n$ with $n\geq 5$, if $|F|\leq 3n-8$ and $Q_n-F$ satisfies the following two conditions:
\begin{enumerate}
  \item [\rm (1)] each vertex has degree at least $2$ in $Q_n-F$;
  \item [\rm (2)] there do not exist a pair of nonadjacent vertices $u$ and $v$ in a $4$-cycle of $Q_n-F$ such that both $u$ and $v$ have degree $2$,
\end{enumerate}
then there is a cycle of every even length from $4$ to $2^n$ in $Q_n-F$.
\end{prop}

Based on Proposition~\ref{prop=5.2}, we have the following lemma.

\begin{lem}\label{lem=5-2}
In a $($faulty$)$ hypercube $Q_n$ with $n\geq 5$, if $|F|\leq 3n-8$ and there exist exactly three vertices of degree $2$ in $Q_n-F$, then $Q_n-F$
has a cycle of every even length from $4$ to $2^n$.
\end{lem}

\f {\bf Proof.} Let $u,v,w$ be the three vertices of degree 2 in $Q_n-F$. Since $|F|\leq 3n-8$ and $Q_n$ is an $n$-regular graph, these three vertices form a faulty 2-path in $Q_n$, and without loss of generality, we assume that it is $(u,v,w)$. Each vertex has degree at least 3 in $Q_n-F$, except for $u,v,w$.

The vertex $v$ has distance 1 with $u$ or $w$. The vertices $u$ and $w$ have distance 2. Since there exists a unique 4-cycle
passing through $u$ and $w$ in $Q_n$, which should contain the faulty 2-path $(u,v,w)$, it follows that there do not exist a 4-cycle containing $u$ and $w$ in $Q_n-F$. Now, $Q_n-F$
satisfies the two conditions in Proposition~\ref{prop=5.2}, and so it has a cycle of every even length from 4 to $2^n$.\hfill\qed

The main result of this section follows directly from the subsequent two lemmas.

\begin{lem}\label{lem=5-1}
In a $($faulty$)$ augmented hypercube $AQ_n$ with $n\geq 5$ and $|F|\leq 4n-8$.
If there exists a vertex $u$ with $d_{AQ_n-F}(u)=2$, then $AQ_n-F$ has a cycle of every even length from $4$ to $2^n$.
\end{lem}

\f {\bf Proof.} Let $F'$ be the subset of edges in $F$ that are not incident to $u$, and let $A=\{u\}\cup F'$. Since $u$ is incident to exactly $2n-3$ faulty edges and the total number of faulty edges satisfies $|F|\leq 4n-8$,
it follows that $|F'|\leq 2n-5$. Consequently, we have $|A|\leq 2n-4$.
Let $(u,x)$ and $(u,y)$ be the two fault-free edges that are incident to $u$. By Proposition~\ref{prop=5.3-1}, there exists a $(u,v)$-path $P[x,y]$ of length $\ell$ ranging from 4 to $2^n-2$ in $AQ_n-A$. Since $E(AQ_n-A)\subset E(AQ_n-F)$,
the path $P[u,v]$ is fault-free. Let $C=(u,x)+P[x,y]+(y,u)$. Then $C$ is a cycle of length $\ell+2$ ranging from 6 to $2^n$ in $AQ_n-F$.

By Proposition~\ref{prop=2.3-2}, there are $2^{n-2}(2n^2+5n-11)$ $4$-cycles in $AQ_n$,
and for each edge $e\in F$, there exist at most
$2n+8$ $4$-cycles going through $e$. Since
$|F|\leq 4n-8$, there are at most
$(2n+8)(4n-8)$ faulty $4$-cycles in $AQ_n$,
and since $2^{n-2}(2n^2+5n-11)>(2n+8)(4n-8)$
when $n\geq5$, there exists at least one 4-cycle in $AQ_n$ that is fault-free.\hfill\qed

The following lemma is a key result of this section, and its proof proceeds in three main steps:
\begin{itemize}
\item  {\bf $Q_n$-isomorphic subgraphs Selection.} Consider the two $Q_n$-isomorphic subgraphs $Q_n^1$ and $Q_n^2$ in $AQ_n$ (see Lemma ~\ref{lem=2}).
\item {\bf Construct fault-free cycles via Proposition~\ref{prop=5.2}.} If $Q_n^1-F$ satisfies the conditions of Proposition~\ref{prop=5.2}, then it has a cycle $C$ of every even length from 4 to $2^n$. Since $Q_n^1-F$ is a subgraph of $AQ_n-F$, the cycle $C$ is automatically a fault-free cycle of $AQ_n$, completing the proof.
\item {\bf Modify $Q_n^1$ via PM-Reciprocity Construction Method.} If $Q_n^1-F$ fails to satisfy Proposition~\ref{prop=5.2}, then apply the PM-Reciprocity Construction Method (see Theorems~\ref{the=3.2-1} and \ref{the=3.2-2}) to modify certain perfect matchings in $Q_n^1$. This constructs a new $Q_n$-isomorphic subgraph $\overline{Q_n^1}$ so that $\overline{Q_n^1}-F$ now meets the requirement of Proposition~\ref{prop=5.2}. Subsequently, applying Proposition~\ref{prop=5.2} to $\overline{Q_n^1}-F$ yields fault-free cycles for $AQ_n$.
\end{itemize}

\begin{lem}\label{lem=5-3}
In a $($faulty$)$ augmented hypercube $AQ_n$ with $n\geq 5$ and $|F|\leq 4n-8$. If each vertex is incident to at least three fault-free edges, then $AQ_n-F$ has a Hamiltonian cycle.
\end{lem}

\f {\bf Proof.} We first choose two $Q_n$-isomorphic subgraphs of $AQ_n$,
whose edge-sets have intersection a perfect matching. Recall that the edge-set of $AQ_n$ can be partitioned to $2n-1$ perfect matchings, i.e., $E(AQ_n)=(\cup_{i=1}^n E_i) \cup (\cup_{j=2}^n E_{\leq j})$. In view of Lemma~\ref{lem=2} and Theorem~\ref{the=3.2-2}, we may assume that $E_1$
is the one with the minimum number of faulty edges in these perfect matchings, and in view of Theorem~\ref{the=3.2-1}, we may further assume that $|E_j\cap F|\leq |E_{\leq j}\cap F|$
for $2\leq j\leq n$.
Since $|F|\leq 4n-8< 2(2n-1)$, we have $|E_1\cap F|\leq 1$. The above assumptions are summarized below:
\begin{equation}\label{eq=5-1}
|E_1\cap F|\leq 1,~|E_1\cap F|\leq |E_j\cap F|\leq |E_{\leq j}\cap F|,~\forall 2\leq j\leq n.
\end{equation}
Now, we consider the $Q_n$-isomorphic subgraphs $Q_n^1$ and $Q_n^2$
in $AQ_n$, that are the edge-induced subgraphs by $E_1\cup \cdots \cup E_n$ and
$E_1\cup E_{\leq 2}\cup \cdots \cup E_{\leq n}$, respectively(see Lemma~\ref{lem=2}). Let $f_i=|F\cap E(Q_n^i)|$ with $i=1$ or $2$. Then $4n-8\geq |F|=f_1+f_2-|E_1\cap F|\geq f_1+f_2-1\geq 2f_1-1$, yielding that $f_1\leq 2n-4$.

To prove this lemma, it suffices to construct a $Q_n$-isomorphic subgraph in $AQ_n$ satisfying the conditions of Proposition~\ref{prop=5.2}.
We consider three cases depending on the minimum degree $d$ of $Q_n^1-F$.

\medskip
\f {\bf Case 1.} $d\geq 2$ .

\begin{figure}
	\centering	\includegraphics[width=0.6\linewidth]{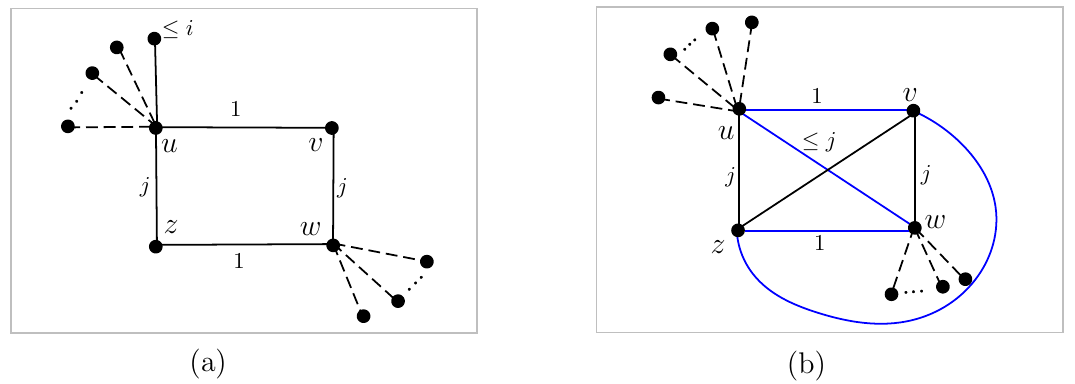}
	\caption{Illustration of  Case 1 of the proof of Lemma~\ref{lem=5-3}}
	\label{fig=1.2}
\end{figure}

Note that $f_1\leq 2n-4\leq 3n-8$ when $n\geq 5$. By Proposition~\ref{prop=5.2},
we may assume that there exists a 4-cycle $C=(u,v,w,z,u)$ in $Q_n^1-F$
such that $d_{Q_n^1-F}(u)=d_{Q_n^1-F}(w)=2$. In this case, $f_1=2n-4$ and all faulty edges in $Q_n^1$ are incident to $u$ or $w$. Each vertex in $Q_n^1-F$ has degree at least $n-1$, except for $u$ and $v$. Given that $|E_1\cap F|\leq 1$,  it is straightforward to observe that $C$ must contains exactly two edges of dimension 1. We may assume that the edges $(u,v)$ and $(z,w)$ are both hypercube edges of dimension 1, and $(u,z)$ and $(v,w)$ are both hypercubes edges of dimension $j$ with $2\leq j\leq n$, see Figure~\ref{fig=1.2} (a). Since all faulty edges in $Q_n^1$ are incident to $u$ or $w$, we have that $|E_i\cap F|=2$ for each $i$ with $2\leq i\neq j\leq n$, and $E_j\cap F=E_1\cap F=\emptyset$.
Since $f_1+f_2-|E_1\cap F|\leq 4n-8$ and $|E_i\cap F|\leq |E_{\leq i}\cap F|$,
we have $f_2=2n-4$, $|E_{\leq i}\cap F|=|E_i\cap F|=2$ and $E_{\leq j}\cap F=\emptyset$.

Now, we begin to replace some perfect matchings in $Q_n^1$ to construct new $Q_n$-isomorphic subgraphs of $AQ_n$.
If $u$ (or $w$) is incident to a
fault-free edge in $E_{\leq i}$ for some $i\neq j$,
then we replace $E_{i}$ in $Q_n^1$ by $E_{\leq i}$, resulting in a subgraph $\overline{Q_n^1}$. By Theorem~\ref{the=3.2-1}, $\overline{Q_n^1}\cong Q_n$. Since $|E_i\cap F|=|E_{\leq i}\cap F|$, we have that $|E(\overline{Q_n^1})\cap F|=|E(Q_n^1)\cap F|=2n-4$. Now, $d_{\overline{Q_n^1}-F}(u)=3$,
$d_{\overline{Q_n^1}-F}(w)\geq 2$, and $d_{\overline{Q_n^1}-F}(x)\geq n-2\geq 3$
for $x\in V(\overline{Q_n^1})\setminus\{u,w\}$. This implies that at most one vertex in $\overline{Q_n^1}-F$ has degree 2, i.e., $\overline{Q_n^1}$ satisfies the condition (2) in Proposition~\ref{prop=5.2}, as required.

If neither $u$ nor $w$ is incident to a
fault-free edge in $E_{\leq i}$ for $i\neq j$,
then all faulty edges in $AQ_n$ are incident to $u$ or $w$.
In particular, $E(AQ_n-\{u,w\})\cap F=\emptyset$.
In $Q_n^2-F$, the vertices $u$ and $w$ are the only vertices having degree 2, and each of them is incident to a hypercube edge of dimension 1 and an augmented edge of dimension $j$. Note that $f_2=2n-4\leq 3n-8$.
If these four edges cannot form a 4-cycle, then $Q_n^2$
satisfies the conditions of Proposition~\ref{prop=5.2}, as required.
Recall that the hypercube edges of dimension 1 incident to $u$ and $w$ are
$(u,v)$ and $(w,z)$, respectively. If these four edges form a 4-cycle, then
$(u,w)$ and $(v,z)$ are augmented edges of dimension $j$ in $Q_n^2-F$, yielding
a fault-free path $(v,u,w,z)$ of length 3 in $AQ_n$, see Figure~\ref{fig=1.2} (b).
By Proposition~\ref{prop=5.3-1}, there exists a $(v,z)$-path $P[v,z]$ in $AQ_n-\{u,w\}$
of length $\ell$ with $3\leq \ell\leq 2^n-3$.
Let $C'=P[v,z]+(v,u,w,z)$. Then $C'$ is a cycle of length $\ell+3$, where $6\leq \ell+3\leq 2^n$. Since $E(AQ_n-\{u,w\})\cap F=\emptyset$, the path $P[v,z]$
is a fault-free path of $AQ_n$, and since $(v,u,w,z)$ is fault-free,
the cycle $C'$ is a fault-free cycle.
Recall that both $Q_n^1$ and $Q_n^2$ have already a fault-free 4-cycle. The proof of Case 1 is complete.

\medskip
\f {\bf Case 2.} $d=1$.

Assume $d_{Q_n^1-F}(u)=1$, i.e., the vertex $u$ is incident to only one fault-free edge in $Q_n^1$, which is a hypercube edge of dimension $i_1$. Since $f_1\leq 2n-4$,
for each vertex $v\neq u$,
\begin{itemize}
\item[\rm (2.1)] either $d_{Q_n^1-F}(v)\geq3$,
or there exists a unique vertex $v$ such that $d_{Q_n^1-F}(v)=2$, $(u,v)\in F$ and all faulty edges in $Q_n^1$ are incident to $u$ or $v$.
\end{itemize}
If $i_1=1$, then $u$ is incident to a faulty edge of dimension $j$ in $Q_n^1$ for each $2\leq j\leq n$, yielding that $1\leq |E_j\cap F|\leq |E_{\leq j}\cap F|$.  If $i_1\neq 1$, then $u$ is incident to a faulty edge of dimension $1$ in $Q_n^1$.
If forces from Eq.~\eqref{eq=5-1} that $|E_{\leq j}\cap F|\geq |E_j\cap F|\geq |E_1\cap F|=1$.
Anyway, we have that $|E_{\leq j}\cap F|\geq |E_j\cap F|\geq 1$ for $2\leq j\leq n$. In $Q_n^2$, the vertex $u$ is incident to at least two fault-free augmented edges. Pick two of such edges and assume that they have dimension $j_1$ and $j_2$ with $2\leq j_1<j_2\leq n$. Without loss of generality, we may assume that $|E_{\leq j_1}\cap F|\leq |E_{\leq j_2}\cap F|$.

Assume $j_1\neq i_1$. Through replacing the $E_{j_1}$ in $Q_n^1$ by $E_{\leq j_1}$, we obtain a $Q_n$-isomorphic subgraph $G_1$ in $AQ_n$ by Theorem~\ref{the=3.2-1}, i.e., the induced subgraph by $E_{1}\cup \cdots \cup E_{j_1-1}\cup E_{\leq j_1}\cup E_{j_1+1}\cup \cdots \cup E_n$. If $|E_{\leq j_1}\cap F|=1$, then $|E_{j_1}\cap F|=1$, and so $|E(G_1)\cap F|=f_1-|E_{j_1}\cap F|+|E_{\leq j}\cap F|=f_1\leq 2n-4\leq 3n-8$ when $n\geq5$.
If $|E_{\leq j_1}\cap F|\geq 2$, then $|E_{\leq j_2}\cap F|\geq2$, and so $|E(G_1)\cap F|=|F|-\sum_{j\neq j_1,j_2}|E_{\leq j_1}\cap F|-|E_{\leq j_2}\cap F|-|E_{j_1}\cap F|\leq 4n-8-(n-3)-2-1\leq 3n-8$.

Assume $j_1=i_1$. It implies that $i_1\geq 2$
and $u$ is incident to a faulty edge of dimension 1. Since $|E_1\cap F|\leq 1$,
it forces that $|E_1\cap F|=1$. Through replacing the $E_{1}$ in $Q_n^1$ by $E_{\leq j_1}$, we obtain a $Q_n$-isomorphic subgraph $G_2$ in $AQ_n$ by Theorem~\ref{the=3.2-2} (Eq.\eqref{eq=5}), i.e, the induced subgraph by $E_{2}\cup \cdots \cup E_{j_1}\cup E_{\leq j_1}\cup E_{j_1+1}\cup \cdots \cup E_n$.
If $|E_{\leq j_1}\cap F|=1$, then $|E(G_2)\cap F|=f_1-|E_{1}\cap F|+|E_{\leq j}\cap F|=f_1-1+1\leq 2n-4\leq 3n-8$ when $n\geq5$.
If $|E_{\leq j_1}\cap F|\geq 2$, then $|E_{\leq j_2}\cap F|\geq2$, and so $|E(G_2)\cap F|=|F|-\sum_{j\neq j_1,j_2}|E_{\leq j_1}\cap F|-|E_{\leq j_2}\cap F|-|E_{1}\cap F|\leq 4n-8-(n-3)-2-1\leq 3n-8$.

Anyway, we can always obtain a new $Q_n$-isomorphic subgraph, say $\overline{Q_n^1}$, by replacing $E_{j_1}$
or $E_1$ in $Q_n^1$ with $E_{\leq j_1}$,
such that $|E(\overline{Q_n^1})\cap F|\leq 3n-8$. Moreover, $d_{\overline{Q_n^1}-F}(u)=2$, and $u$ is incident to a hypercube edge of dimension $i_1$ and an augmented edge of dimension $j_1$.
Among the edges incident to $v$ with $v\neq u$, by the construction of $\overline{Q_n^1}$,
only the hypercube edge of dimension $j_1$ or 1 is changed to be an augmented edge of dimension $j_1$, and the other edges keep unchanged.

{\bf Subcase 2.1} Each vertex has degree at least 2 in $\overline{Q_n^1}-F$.

Since $|E(\overline{Q_n^1})\cap F|\leq 3n-8$, there are at most three vertices of degree 2 in $\overline{Q_n^1}-F$, and by Lemma~\ref{lem=5-2} and Proposition~\ref{prop=5.2}, we may assume that
there are exactly two vertices of degree 2 and they are not adjacent. One of them is $u$ and the other one is denoted to be $v$.
It follows from (2.1) that the augmented edge of dimension $j_1$ incident to $v$ is faulty, i.e., $v$ is incident to two hypercube edges in $\overline{Q_n^1}-F$.
Recall that $u$ is incident to a hypercube edge and an augmented edge of dimension $j_1$ in $\overline{Q_n^1}-F$. Suppose that there exists a 4-cycle $C$ in $\overline{Q_n^1}-F$ containing $u$ and $v$. Then $C$ has only one augmented edge which is of dimension $j_1$. By Proposition~\ref{prop=2.3-3}~(1), $j_1=3$ and $C$ contains a hypercube edge of dimension $1$
and a hypercube edge of dimension $j_1=3$, which is impossible because one of $E_1$ and $E_{j_1}$ does not belong to $E(\overline{Q_n^1})$.
Therefore, $\overline{Q_n^1}-F$ satisfies that second condition of Proposition~\ref{prop=5.2}, as required.

{\bf Subcase 2.2} There exists a vertex $v$ having degree at most 1 in $\overline{Q_n^1}-F$.

Recall that from $Q_n^1$ to $\overline{Q_n^1}$,
only one hypercube edge of dimension $j_1$ or 1 incident to $v$ is changed.
It follows that $v$ has degree at most 2 in $Q_n^1-F$, and from (2.1) that such vertex $v$ exists uniquely with $(u,v)\in F\cap E(Q_n^1)$ and all faulty edges in $Q_n^1$ are incident to $u$ or $v$. The vertex $x$ is incident to at least $n-1\geq4$ fault-free edges in $Q_n^1$ when $x\neq u,v$.

Now, we revisit the search for a $Q_n$-isomorphic subgraph within $AQ_n$, employing the same approach of modifying a perfect matching in $Q_n^1$. We have that
$u$ is incident to at least 3 fault-free edges
in $AQ_n$, three of which are a hypercube edge of dimension $i_1$ and two augmented edges of dimensions $j_1$ and $j_2$. The vertex $v$ is also incident to at least 3 fault-free edges
in $AQ_n$, two of which are hypercube edges that belong to $Q_n^1$. Assume that
they have dimensions $i_2$ and $i_3$ with $i_2<i_3$, and let $(u,v)$ is of dimension $i_4$. Since $|E_{1}\cap F|\leq 1$,
at least one of $i_1$ and $i_2$ is equal to 1.
Recall that $|E_j\cap F|\leq |E_{\leq j}\cap F|$ for $2\leq j\leq n$.

\begin{figure}
	\centering
	\includegraphics[width=0.6\linewidth]{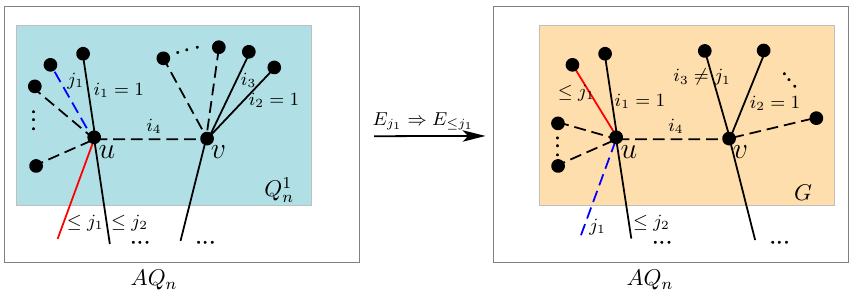}
	\caption{Illustration of Subcase 2.2 with $i_1=i_2=1$ in the proof of Lemma~\ref{lem=5-1}}
	\label{fig=2.2}
\end{figure}

Assume $i_1=i_2=1$. Then $|E_{i_3}\cap F|=|E_{i_4}\cap F|=1$ and $|E_{j}\cap F|=2$ for $2\leq j\neq i_3,i_4\leq n$. Without loss of generality, we may assume that $j_1\neq i_3$. Exchanging $E_{j_1}$ in $Q_n^1$ with $E_{\leq j_1}$, and obtaining a $Q_n$-isomorphic subgraph of $AQ_n$ by Theorem~\ref{the=3.2-1}, say $G$(see Figure~\ref{fig=2.2}).
Consequently, $|E(G)\cap F|=|F|-\sum_{j\neq j_1,i_3,i_4} |E_{\leq j}\cap F|-|E_{\leq i_3}\cap F|-|E_{\leq i_4}\cap F|-|E_{j_1}\cap F|\leq 4n-8-2(n-4)-1-1-1=2n-3\leq 3n-8$. (Note that $j_1$ may be equal $i_4$.) Moreover,
$d_{G-F}(u)=2$, $d_{G-F}(v)\geq2$ and $d_{G-F}(x)\geq 3$ when $x\neq u,v$. In $G-F$, the vertex $u$ is incident to an augmented edge of dimension $j_1$ and a hypercube edge of dimension $i_1=1$; the vertex $v$ is incident to two hypercube edges of dimensions 1 and $i_3$.
If $d_{G-F}(v)\geq3$, then $G-F$ has only one vertex of degree 2, and so satisfies the conditions of Proposition~\ref{prop=5.2}, as required. If $d_{G-F}(v)=2$, then by a similar argument as Subcase 3.1, there does not exist a 4-cycle containing $u$ and $v$. Hence $G$ satisfies the conditions of Proposition~\ref{prop=5.2}, as required.

\begin{figure}
	\centering
	\includegraphics[width=0.6\linewidth]{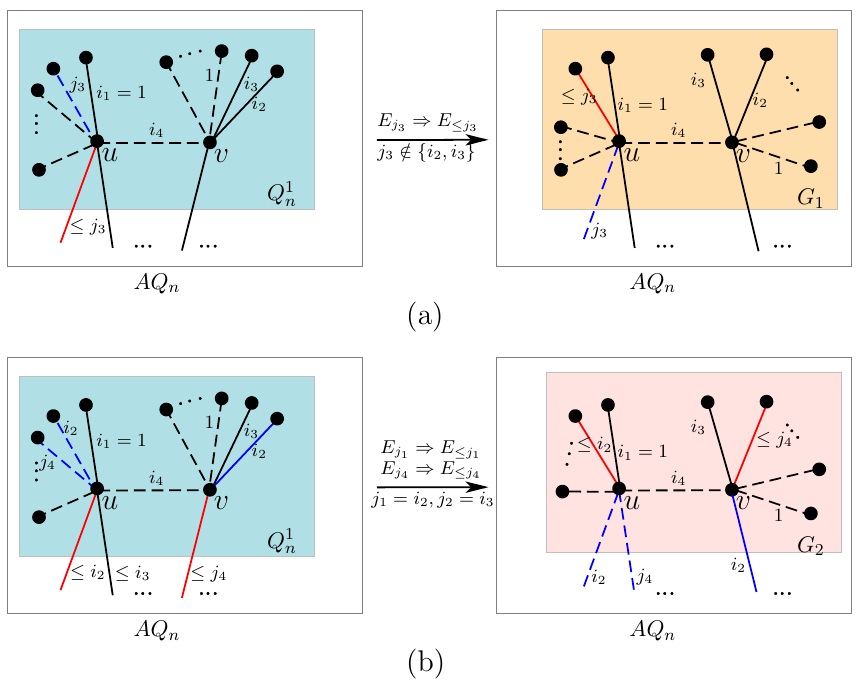}
	\caption{Illustration of Subcase 2.2 with $i_1=1$ and $i_2\neq 1$ in the proof of Lemma~\ref{lem=5-1}}
	\label{fig=2.4}
\end{figure}

Assume $i_1=1$ and $i_2\neq1$.
Then $|E_{i_1}\cap F|=|E_{i_2}\cap F|=|E_{i_3}\cap F|=|E_{i_4}\cap F|=1$ and $|E_{j}\cap F|=2$ when $j\notin\{i_1,i_2,i_3,i_4\}$.
If there is a $j_3$ such that $j_3\notin\{i_2,i_3\}$ and $u$ is incident to a fault-free augmented edge of dimension $j_3$ (see Figure~\ref{fig=2.4} (a)),
then replacing $E_{j_3}$ in $Q_n^1$ with $E_{\leq j_3}$ and resulting in a $Q_n$-isomorphic subgraph $G_1$ of $AQ_n$ by Theorem~\ref{the=3.2-1}.
Consequently,
$|E(G_1)\cap F|=|F|-\sum_{j\neq i_2,i_3,i_4,j_3} |E_{\leq j}\cap F|-|E_{\leq i_2}\cap F|-|E_{\leq i_3}\cap F|-|E_{\leq i_4}\cap F|-|E_{j_3}\cap F|\leq 4n-8-2(n-5)-1-1-1-2=2n-3\leq 3n-8$ when $j_3\neq i_4$ and $|E(G_1)\cap F|\leq 4n-8-2(n-4)-1-1-1=2n-3$ when $j_3=i_4$. Moreover,
$d_{G_1-F}(u)=2$, $d_{G_1-F}(v)\geq2$ and $d_{G_1-F}(x)\geq 3$ when $x\neq u,v$. In $G_1-F$, the vertex $u$ is incident to a hypercube edge of dimension $i_1=1$ and an augmented edge of dimension $j_3$, while $v$ is incident to two hypercube edges of dimensions $i_2$ and $i_3$ with $j_3\notin\{i_2,i_3\}$. By a similar argument as above paragraph,
it follows from Proposition~\ref{prop=2.3-3}~(1) that there does not exist a 4-cycle containing $u$ and $v$. Hence $G_1$ satisfies the conditions of Proposition~\ref{prop=5.2}, as required.

If there does not exist a $j_3$ such that $j_3\notin\{i_2,i_3\}$ and $u$ is incident to a fault-free augmented edge of dimension $j_3$, then $j_1=i_2$ and $j_2=i_3$ (see Figure~\ref{fig=2.4}~(b)).
Since $d_{AQ_n-F}(v)\geq3$, $v$ is incident a fault-free augmented edge of some dimension $j_4$ ($j_4$ may be equal $j_1$ or $j_2$).
Replacing $E_{j_1}$ and $E_{j_4}$ by
$E_{\leq j_1}$ and $E_{\leq j_4}$,
and by Theorem~\ref{the=3.2-1}, it results in a $Q_n$-isomorphic subgraph $G_2$ of $AQ_n$. Consequently, $|E(G_2)\cap F|=|F|-\sum_{j\neq j_1,j_2,j_4}|E_{\leq j}\cap F|-|E_{j_1}\cap F|-|E_{j_4}\cap F|-|E_{\leq j_2}\cap F|\leq
4n-8-2(n-4)-1-1-1=2n-3$ when $j_4\notin\{j_1,j_2\}$, and $|E(G_2)\cap F|\leq |F|-\sum_{j\neq j_1,j_2}|E_{\leq j}\cap F|-|E_{j_1}\cap F|-|E_{j_2}\cap F|\leq
4n-8-2(n-3)-1-1=2n-4$ when $j_4\in\{j_1,j_2\}$. Moreover,
$d_{G_2-F}(u)\geq 2$, $d_{G_2-F}(v)=2$ and $d_{G_2-F}(x)\geq 3$ when $x\neq u,v$. In $G_2-F$, the vertex $u$ is incident to a hypercube edge of dimension $1$ and an augmented edge of dimension $j_1$, while $v$ is incident to a hypercube edge of dimension $i_3=j_2$ and an augmented edge of dimension $j_4$. Suppose that $d_{G_2-F}(u)=2$ and there is a 4-cycle $C$ containing $u$ and $v$.
Then $C$ has exactly two hypercube edges, one of which is is of dimension 1 and the other one is of dimension $i_3\neq 1$, which is contradict to
Proposition~\ref{prop=2.3-2}~(2).
Finally, $G_2$ satisfies the conditions of Proposition~\ref{prop=5.2}, as required.

The case $i_1\neq 1$ and $i_2=1$ is similar as $i_1=1$ and $i_2\neq1$.

\medskip

\f {\bf Case 3.} $d=0$.

There exists a vertex in $Q_n^1-F$ of degree 0. Since $f_1\leq 2n-4$ and
$Q_n^1$ has valency $n$, such vertex in $Q_n^1-F$ is unique, say $u$. All hypercube edges incident to $u$ are faulty, implying that $|E_{i}\cap F|\geq 1$ for $1\leq i\leq n$. In particular, $u$ is incident a faulty edge of dimension 1.
Since $|E_1\cap F|\leq 1$, for each vertex $v\neq u$,
\begin{itemize}
\item [\rm (3.1)] either $(u,v)\in E_1\cap F$ or $v$ is incident to a fault-free edge of dimension 1.
\end{itemize}
Note that $f_1\leq 2n-4$.
Each vertex $v\in V(Q_n^1-\{u\})$ is incident to at most $n-3$ faulty edges in $Q_n^1$, and specially, if $v$ is incident to exactly $n-3$ faulty edges in $Q_n^1$,
then $(u,v)\in F\cap E(Q_n^1)$ and such $v$ exists uniquely. In other words, in $Q_n^1-\{u\}$ either \begin{itemize}
  \item [\rm (3.2)] each vertex $v$ is incident to at least $4$ fault-free edges in $Q_n^1$; or
  \item [\rm (3.3)] there exists a unique vertex $v$ such that it is incident to exactly $3$ fault-free edges in $Q_n^1$ and  $(u,v)\in F\cap E(Q_n^1)$, and each of the other vertices is incident to at least $4$ fault-free edges in $Q_n^1$. All faulty edges in $Q_n^1$ are incident to $u$ or $v$.
\end{itemize}

We now consider the fault-free edges in $AQ_n$ that are incident to $u$. There are at least three such edges, and all of them are augmented edges. Pick three of them and assume that they have dimensions $j_1$, $j_2$ and $j_3$ with $2\leq j_1,j_2,j_3\leq n$ and
$|E_{\leq j_1}\cap F|\leq |E_{\leq j_2}\cap F|\leq |E_{\leq j_3}\cap F|$.
By Theorem~\ref{the=3.2-2}, there exist two $Q_n$-isomorphic subgraphs of $AQ_n$, say $\overline{Q_n^1}$ and $\overline{Q_n^2}$,
such that $E(\overline{Q_n^1})\cap E(\overline{Q_n^2})=E_1$ and $E(\overline{Q_n^1})\cup E(\overline{Q_n^2})=E(AQ_n)$. Here, $\overline{Q_n^1}$
is derived from $Q_n^1$ by replacing the perfect matchings $E_{j_1}$ and $E_{j_2}$ with
$E_{\leq j_1}$ and $E_{\leq j_2}$, while $\overline{Q_n^2}$ is obtained from $Q_n^2$ by substituting $E_{\leq j_1}$ and $E_{\leq j_2}$
with $E_{j_1}$ and $E_{j_2}$. Consequently, $u$ is incident to exactly two fault-free edges in $\overline{Q_n^1}$. Let $\bar{f_i}=|E(\overline{Q_n^i})\cap F|$ for
$i=1$ or 2.

We claim that $\bar{f_1}\leq 3n-8$.
If $|E_{\leq j_1}\cap F|+|E_{\leq j_2}\cap F|\leq2$, then $|E_{\leq j_1}\cap F|+|E_{\leq j_2}\cap F|=|E_{j_1}\cap F|+|E_{j_2}\cap F|=2$
because of $1\leq |E_j\cap F|\leq |E_{\leq j}\cap F|$ for $2\leq j\leq n$.
It follows that
$\bar{f_1}=f_1+|E_{\leq j_1}\cap F|+|E_{\leq j_2}\cap F|-|E_{j_1}\cap F|-|E_{j_2}\cap F|=f_1\leq 2n-4\leq 3n-8$.
If $|E_{\leq j_1}\cap F|+|E_{\leq j_2}\cap F|>2$, then $|E_{\leq j_3}\cap F|\geq2$
because of $|E_{\leq j_1}\cap F|\leq |E_{\leq j_2}\cap F|\leq |E_{\leq j_3}\cap F|$.
Therefore, $\bar{f_1}=|F|-\sum_{j\neq j_1,j_2,j_3} |E_{\leq j}\cap F|-|E_{\leq j_3}\cap F|-|E_{j_1}\cap F|-|E_{j_2}\cap F|\leq 4n-8-(n-4)-2-1-1=3n-8$. As claimed.

Now, $u$ is incident to exactly two fault-free edges in $\overline{Q_n^1}$, which are augmented edges of dimensions $j_1$ and $j_2$ respectively. For each vertex $v$ with $v\neq u$, only two edges incident to $v$ in $Q_n^1$ are changed in the construction of $\overline{Q_n^1}$. It follows from (3.2) and (3.3) that either each vertex in $\overline{Q_n^1}-\{u\}$ is incident to at least 2 fault-free edges in $\overline{Q_n^1}$; or there exists a unique vertex $v$ in $\overline{Q_n^1}$ such that $v$ is incident to exactly 1 fault-free edges and $(u,v)\in F\cap E(Q_n^1)$, and the other vertices are all incident to more than 2 fault-free edges in $\overline{Q_n^1}$.

For the former case, $\overline{Q_n^1}$
satisfies the condition (1) in Proposition~\ref{prop=5.2}. Since $\bar{f_1}\leq 3n-8$, there exist at most three vertices of degree 2 in $\overline{Q_n^1}-F$, one of which is $u$. If there exists three such vertices, then the lemma holds by Lemma~\ref{lem=5-2}. If there exist exactly two vertices of degree 2 in $\overline{Q_n^1}-F$ which are adjacent, or there exists only one vertex of degree 2 in $\overline{Q_n^1}-F$, then $\overline{Q_n^1}$ satisfies the condition (2) in Proposition~\ref{prop=5.2}.
Assume that there exists exactly two vertices of degree 2 in $\overline{Q_n^1}-F$, say $u$ and $v$, such that $u$ is not adjacent to $v$. The two edges incident to $u$ in $\overline{Q_n^1}-F$ are augmented edges, but one edge incident to $v$ in $\overline{Q_n^1}-F$ belongs to $E_1$ (see (2.1)). It follows from Proposition~\ref{prop=2.3-3}~(3) and (4) that there does not exist a 4-cycle containing $u$ and $v$ in $\overline{Q_n^1}-F$, i.e., $\overline{Q_n^1}$
satisfies the condition (2) in Proposition~\ref{prop=5.2}, as required.

\begin{figure}
	\centering
	\includegraphics[width=0.9\linewidth]{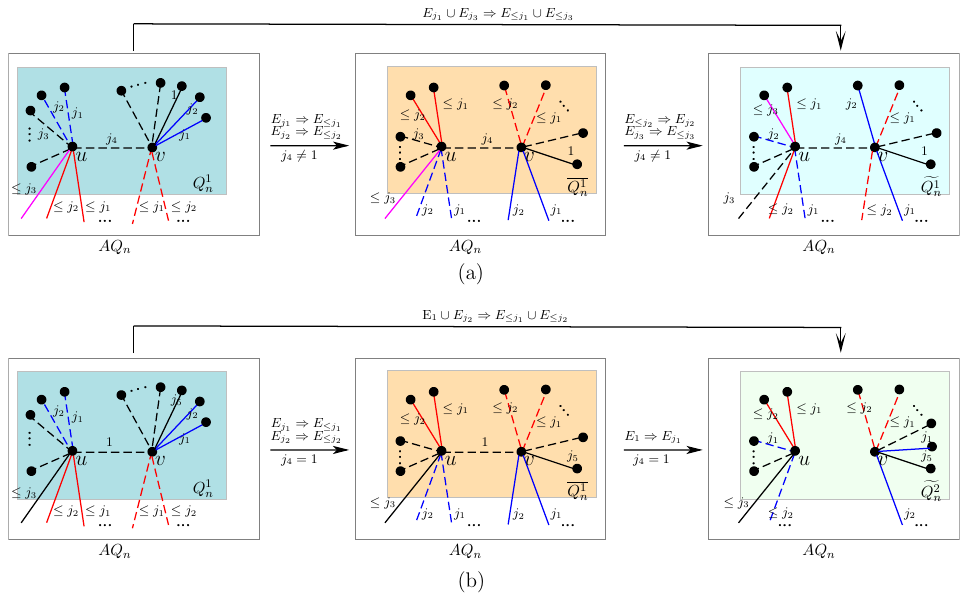}
	\caption{Illustration of Case 3 in the proof of Lemma~\ref{lem=5-1}}
	\label{fig=2.5}
\end{figure}

For the latter case, there exists a unique vertex $v$ in $\overline{Q_n^1}$ such that $v$ is incident to exactly 1 fault-free edges and $(u,v)\in F\cap E(Q_n^1)$.
Recall that in the construction of $\overline{Q_n^1}$, the $E_{j_1}\cup E_{j_2}$
in $Q_n^1$ are replaced by $E_{\leq j_1}\cup E_{\leq j_2}$. In $Q_n^1$, $v$ is incident to 3 faulty edges (see (2.3)), and now in $\overline{Q_n^1}$ it is incident to 1 faulty edges. It implies that the edges incident to $v$ in $E_{j_1}\cup E_{j_2}$ are fault-free, while in $E_{\leq j_1}\cup E_{\leq j_2}$ are faulty (see Figure~\ref{fig=2.5}~(a)).
Assume that $(u,v)$ is of dimension $j_4$.
If $j_4\neq 1$, then by (3.1) the three fault-free edges incident to $v$ in $Q_n^1$ are hypercube edges of dimension 1, $j_1$ and $j_2$, respectively.
Since all faulty edges in $Q_n^1$ are incident to $u$ or $v$ (see (3.3)), we have that
$|E_{j_1}\cap F|=|E_{j_2}\cap F|=|E_{j_4}\cap F|=1$ and
$|E_j\cap F|=2$ with $2\leq j\neq j_1,j_2,j_4\leq n$.
We restart modifying the subgraph $Q_n^1$.
Let $\widetilde{Q_n^1}$ be the subgraph of $AQ_n$ induced by
$$E_1\cup \cdots \cup E_{j_1-1}\cup E_{\leq j_1}\cup E_{j_1+1}\cup \cdots \cup E_{j_3-1}\cup E_{\leq j_3}\cup E_{j_3+1}\cup \cdots \cup E_n.$$
The $\widetilde{Q_n^1}$ is derived from $Q_n^1$ (or $\overline{Q_n^1}$, resp.) by replacing the perfect matching $E_{j_1}\cup E_{j_3}$ (or $E_{\leq j_2}\cup E_{j_3}$, resp.) with $E_{\leq j_1}\cup E_{\leq j_3}$ (or $E_{j_2}\cup E_{\leq j_3}$, resp.), see Figure~\ref{fig=2.5}~(a). By Theorem~\ref{the=3.2-1}, $\widetilde{Q_n^1}\cong Q_n$. Let $\widetilde{f_1}=|E(\widetilde{Q_n^1})\cap F|$.
Then $\widetilde{f_1}=|F|-\sum_{j\neq j_1,j_2,j_3,j_4}|E_{\leq j}\cap F|-|E_{\leq j_2}\cap F|-|E_{\leq j_4}\cap F|-|E_{j_1}\cap F|-|E_{j_3}\cap F|\leq 4n-8-2(n-5)-1-1-1-2=2n-3\leq 3n-8$.
After modifying, the vertex $u$ is incident to exactly two edges
in $\widetilde{Q_n^1}-F$ which are augmented edges of dimension $j_1$ and $j_3$; $v$ is incident to a hypercube edge of dimension $j_2$ and
a hypercube edge of dimension 1. For each vertex $x$ with $x\neq u$ or $v$, it has degree at least 4 in $Q_n^1-F$ by (3.3), and so it has degree at least 2 in $\widetilde{Q_n^1}-F$.
It can be checked by a similar argument as in above paragraph that $\widetilde{Q_n^1}$ satisfies the condition (2) in Proposition~\ref{prop=5.2}, as required.

If $j_4=1$, then $v$ is incident to three hypercube edges in $Q_n^1$ of dimensions $j_1$, $j_2$ and $j_5$ for some $j_5\neq 2$,
yielding that $|E_{j_1}\cap F|=|E_{j_2}\cap F|=|E_{j_5}\cap F|=1$ and
$|E_j\cap F|=2$ with $2\leq j\neq j_1,j_2,j_5\leq n$.
Let $\widetilde{Q_n^2}$ be the subgraph induced by
$$E_2\cup \cdots \cup E_{j_1}\cup E_{\leq j_1}\cup E_{j_1+1}\cup \cdots \cup E_{j_2-1}\cup E_{\leq j_2}\cup E_{j_2+1}\cup \cdots \cup E_n.$$
The $\widetilde{Q_n^2}$ is derived from $Q_n^1$ (or $\overline{Q_n^1}$, resp.) by replacing the perfect matching $E_{1}\cup E_{j_2}$ (or $E_{1}$, resp.) with $E_{\leq j_1}\cup E_{\leq j_2}$ (or $E_{j_1}$, resp.), see Figure~\ref{fig=2.5}~(b). By Lemma~\ref{lem=3.2-2}, $\widetilde{Q_n^2}\cong Q_n$. Let $\widetilde{f_2}=|E(\widetilde{Q_n^2})\cap F|$.
Then $\widetilde{f_2}=|F|-|E_1\cap F|-\sum_{j\neq j_1,j_2,j_5}|E_{\leq j}\cap F|-|E_{j_2}\cap F|-|E_{\leq j_5}\cap F|\leq 4n-8-1-2(n-4)-1-1=2n-3\leq 3n-8$.
The vertex $u$ is incident to two augmented edges
in $\widetilde{Q_n^1}-F$ of dimensions $j_1$ and $j_2$ with $j_1<j_2$, while $v$ is incident to two hypercube edges of dimensions $j_1$ and $j_5$. For each vertex $x$ with $x\neq u$ or $v$, it has degree at least 2 in $\widetilde{Q_n^1}-F$.
By a similar argument as in above paragraph that $\widetilde{Q_n^1}-F$ satisfies the two conditions in Proposition~\ref{prop=5.2}, as required.

The proof is complete.
\hfill\qed

Combined with Lemmas~\ref{lem=5-1} and \ref{lem=5-3}, we have our main result as follows.

\begin{theorem}\label{theo=5.1}
In a $($faulty$)$ augmented hypercube $AQ_n$ with $n\geq 5$,  if $|F|\leq 4n-8$ and
each vertex is incident to at least two fault-free edges, then $AQ_n-F$ has a cycle of every even length from $4$ to $2^n$. In particular, $AQ_n-F$ has a Hamiltonian cycle.
\end{theorem}

\section{Conclusion}

The augmented cube $AQ_n$, as a variant of the hypercube $Q_n$, represents an important interconnection network topology with significant applications in parallel and distributed computing. Due to its advantageous structural properties, $AQ_n$ proves particularly suitable for parallel algorithm design, fault-tolerant routing, broadcasting, and related applications~\cite{CS}. In this paper, we conduct a systematic investigation of the structural properties of augmented cubes $AQ_n$, with special emphasis on their spanning subgraphs.

By leveraging the Cayley graph properties of $AQ_n$, we give a method for constructing $Q_n$-isomorphic spanning subgraphs. Our approach specifically focuses on constructing pairs of $Q_n$-isomorphic subgraphs with minimum number of common edges. The construction process involves: (1) partitioning the edge set of $AQ_n$ into $2n-1$ perfect matchings, and (2) strategically selecting perfect matchings through minimal generating sets of the elementary abelian 2-group $\mz_2^n$. We introduce the PM-reciprocity construction method, which enables the generation of additional $Q_n$-isomorphic subgraphs through systematic modifications of existing perfect matchings. Based on these methods, we establish a lower bound for the number of $Q_n$-isomorphic subgraphs in $AQ_n$, thereby making partial progress toward resolving an open problem posed by Dong and Wang~\cite{DW}.

As an application of our method, we establish that $AQ_n$ contains $n-1$ or $n-2$ edge-disjoint Hamiltonian cycles when $n$ is odd or even, respectively. Notably, when $n$ is odd, our result implies that $AQ_n$ admits a Hamiltonian decomposition, thereby confirming the odd case of a conjecture proposed by Hung~\cite{Hung3}. The even case remains open and appears to require substantially different techniques for its resolution.

We also discuss the edge fault-tolerant Hamiltonicity of $AQ_n$. In the end of this paper, we prove that under the conditional fault model, $AQ_n-F$ with $n\geq5$ has a cycle of every even length from 4 to $2^n$ even if $|F|\leq 4n-8$. Notably, this implies that $AQ_n-F$ remains Hamiltonian. This not only provides an alternative proof for the well-known fault-tolerant Hamiltonicity of $AQ_n$ established by Hsieh and Cian~\cite{HC}, but also extends their work by demonstrating the fault-tolerant bipancyclicity of $AQ_n$. Our proof leverages the existence of multiple $Q_n$-isomorphic subgraphs within $AQ_n$, particularly employing the PM-reciprocity construction method. This approach provides a new framework for studying fault-tolerant cycle embeddings in augmented cubes.
This work naturally leads to several important open questions, including extending these results to fault-tolerant pancyclicity (embedding cycles of arbitrary lengths).  We note that our current proof relies on constructing spanning subgraphs isomorphic to the hypercube, which are bipartite.
To study the embedding of odd cycles in $AQ_n$, future work may need to consider non-bipartite spanning subgraphs, such as folded hypercube $FQ_n$. In particular, identifying subgraphs of $AQ_n$ that are isomorphic to $FQ_n$ may deserve further exploration in subsequent work.

\section*{Acknowledgments}
This work was supported by the National Natural Science Foundation of China (Nos. 12101070, 12471321, 12331013, 12271024).

\end{document}